\documentstyle[amssymb,righttag]{amsart}

\advance\textheight by \topskip

\newtheorem{teo}{\sc Theorem}[section]

\newtheorem{cor}{\sc Corollary}[section]
\newtheorem{lemma}{\sc Lemma}[section]
\newtheorem{pro}{\sc Proposition}[section]

\theoremstyle{definition}
\newtheorem{defi}{\sc Definition}[section]

\theoremstyle{remark}
\newtheorem{rem}[teo]{\sc Remark}

\newtheorem{notation}{\sc Notation}

\newcommand{\bte}{\begin{teo}}
\newcommand{\ete}{\end{teo}}
\newcommand{\bc}{\begin{cor}}
\newcommand{\ec}{\end{cor}}
\newcommand{\bp}{\begin{pro}}
\newcommand{\ep}{\end{pro}}
\newcommand{\bl}{\begin{lemma}}
\newcommand{\el}{\end{lemma}}
\newcommand{\bd}{\begin{defi}}
\newcommand{\ed}{\end{defi}}
\newcommand{\bno}{\begin{notation}}
\newcommand{\eno}{\end{notation}}
\newcommand{\fp}{\hfill $\Box$}

\newcommand{\bca}{\begin{cases}}
\newcommand{\eca}{\end{cases}}
\newcommand{\la}{\langle}
\newcommand{\ra}{\rangle}

\newcommand{\bq}{\begin{equation}}
\newcommand{\eq}{\end{equation}}
\newcommand{\btabu}{\begin{table}}
\newcommand{\etabu}{\end{table}}
\newcommand{\bt}{\begin{tabular}}
\newcommand{\et}{\end{tabular}}
\newcommand{\ba}{\begin{array}}
\newcommand{\ea}{\end{array}}
\newcommand{\br}{\begin{eqnarray}}
\newcommand{\er}{\end{eqnarray}}
\newcommand{\brn}{\begin{eqnarray*}}
\newcommand{\ern}{\end{eqnarray*}}
\newcommand{\benu}{\begin{enumerate}}
\newcommand{\eenu}{\end{enumerate}}
\newcommand{\bite}{\begin{itemize}}
\newcommand{\eite}{\end{itemize}}

\newcommand{\supp}{\operatorname{supp }}

\title[Nikishin systems are perfect] {Nikishin systems are perfect. Case of unbounded and touching supports}

\thanks{Both authors were
supported by MICINN of Spain under grants MTM2009-12740-C03-01  and Acciones Integradas Portugal. GL acknowledges that his work was completed while visiting Vanderbilt University on sabbatical from Universidad Carlos III de Madrid. UF acknowledges support from SFRH/BPD/62947/2009, Fundacao para a Ci\^{e}ncia e a Tecnolog\'{i}a of Portugal, PT2009-0031.}

\author[Fidalgo]{U. Fidalgo Prieto}

\address[Fidalgo]{Departamento de Matem\'{a}tica, Universidade de Aveiro, Campus Universit\'{a}rio de Santiago, 3810-193 Aveiro, Portugal}

\email[Fidalgo]{ufidalgo@@math.uc3m.es}

\author[L\'{o}pez]{G. L\'opez Lagomasino}
\address[L\'{o}pez]{Departamento de Matem\'aticas,
Universidad Carlos III de Madrid, c/ Universidad 30, 28911 Legan\'es, Spain.}

\email[L\'{o}pez]{lago@@math.uc3m.es}

\begin{document}

\maketitle

\begin{abstract}
K. Mahler introduced the concept of perfect systems in the theory of simultaneous Hermite-Pad\'{e} approximation of analytic functions.  Recently, we proved that Nikishin systems, generated by measures with bounded support and non-intersecting consecutive supports contained on the real line, are perfect. Here, we prove that they are also perfect when the supports of the generating measures are unbounded or touch at one point. As an application, we give a version of Stieltjes' theorem in the context of simultaneous Hermite-Pad\'e approximation.
\end{abstract}

\vspace{1cm}


{\it A.M.S. Subject Classification.} Primary: 30E10, 42C05;
Secondary: 41A20.

\section{Introduction} \label{sec:HP}

This paper complements \cite{LF4} where we solved a long standing problem proving that Nikishin systems (generated by measures whose supports are bounded and consecutive supports do not intersect) are perfect. Here, we consider Nikishin systems whose generating measures may have unbounded support or consecutive supports touch at one point. A detailed account on the history of the problem is contained in the introduction of \cite{LF4}  so we will go directly to the subject matter.

\subsection{Perfect systems}
The concept of perfect systems was introduced and developed by  K. Mahler in lectures delivered at the University of Groningen in 1934-35 which were published much later in \cite{Mah} (see also \cite{Coa} and \cite{Jag}). This notion plays a central role in the general theory of simultaneous approximation of systems of analytic functions. Systems of exponential functions are perfect and the corresponding properties were used by Ch. Hermite in \cite{Her} to prove the transcendence of $e$. This may have inspired Mahler to introduce the general concept and study its properties.

Mahler's general approach to the simultaneous approximation of finite systems of analytic functions may be reformulated in the following terms.

Let ${\bf f} = (f_0,\ldots,f_m)$ be a  system of formal power expansions at $\infty$ of the form
\[ f_j(z) = \sum_{n=0}^{\infty} \frac{f_{j,n}}{z^n}, \qquad j=0,\ldots,m.
\]
Fix a non-zero multi-index ${\bf n} = (n_0,\ldots,n_m) \in {\mathbb{Z}}_+^{m+1}, |{\bf n}| = n_0+\ldots,n_m$. There exist polynomials $a_{{\bf n},0}, \ldots, a_{{\bf n},m}$, not all identically equal to zero, such that
\begin{itemize}
\item[i)] $\deg a_{{\bf n},j} \leq n_j -1, j=0,\ldots,m$ $(\deg a_{{\bf n},j} \leq -1$ means that $a_{{\bf n},j} \equiv 0$),
\item[ii)] $\sum_{j=0}^m a_{{\bf n},j}(z) f_j(z) - b_{\bf n}(z) = \frac{A_{ \bf n}}{z^{|{\bf n}|}} + \cdots$,
\end{itemize}
for some polynomial $b_{\bf n}$. Analogously, there exists  a polynomial $Q_{\bf n}$, not identically equal to zero, such that
\begin{itemize}
\item[i)] $\deg Q_{{\bf n}} \leq |{\bf n}|,$
\item[ii)] $Q_{\bf n}(z)f_j(z) - P_{{\bf n},j}(z)= \frac{A_{{\bf n},j}}{z^{n_j+1}} + \cdots, \,\, j=0,\ldots,m,$
\end{itemize}
for some polynomials $P_{{\bf n},j}, j=0,\ldots,m.$

The right hand of ii) must be understood as a formal expansion in decreasing powers of $z$ obtained after carrying out arithmetically the operations of the left hand. Certainly, if the left hand is analytic at $\infty$ the right hand will be convergent in some neighborhood of $\infty$ and equality is in the usual sense in that neighborhood.

The polynomials $b_{\bf n}$ and $P_{{\bf n},j},j=0,\ldots,m,$ are uniquely determined from ii) once their partners $a_{{\bf n},j}, j=0,\ldots,m$, and  $Q_{\bf n}$  are found. The two constructions are called type I and type II polynomials (approximants) of the system $(f_0,\ldots,f_m)$, respectively.  When $m=0$ both definitions reduce to the well-known Pad\'{e} approximation in its linear presentation.

In applications (number theory, convergence of simultaneous rational approximation, asymptotic properties of type I and type II polynomials, non-intersecting brownian motions, and random matrix theory) it is important that the polynomials appearing in the construction have no defect; that is, that they have full degree. For one, this guarantees uniqueness up to a constant factor. Here, the following concept steps in.

\begin{defi}\label{def1}
A multi-index ${\bf n} =(n_0,\ldots,n_m)$ is said to be {\bf normal} for the system ${\bf f}$ for type I approximation (respectively, for type II,) if $\deg a_{{\bf n},j} = n_j -1, j=0,\ldots,m$ (respectively, $\deg Q_{\bf n} = |{\bf n}|$). A system of functions ${\bf f}$ is said to be {\bf perfect} if all multi-indices are normal.
\end{defi}

It is easy to see that normality implies that all solutions are collinear.

\subsection{Nikishin systems.} Let $s$ be a finite Borel measure with constant sign whose support consists of infinitely many points contained in the real line. By $\Delta = \mbox{Co}(\supp s)$ we denote the smallest interval which contains $\supp{s},$ the support  of $s$. We denote this class of measures by ${\mathcal{M}}(\Delta)$. Let
\[ \widehat{s}(z) = \int \frac{ds(x)}{z-x}
\]
denote the Cauchy transform of $s$. Obviously, $\widehat{s} \in {\mathcal{H}}(\overline{\mathbb{C}}\setminus \Delta);$  that is, it is analytic in $\overline{\mathbb{C}}\setminus \Delta.$

Assume that all the moments of $s$ are finite; that is,
\[ c_{\nu} = \int x^{\nu} ds(x) \in {\mathbb{R}}, \qquad \nu \in {\mathbb{Z}}_+.
\]
When $s$ is finite and $\supp s$ is bounded this is automatically fulfilled.
If we construct the Pad\'e approximation to $\widehat{s}$, for some index $n \in {\mathbb{Z}}_+$,  $Q_n$ turns out to be the $n$-th orthogonal polynomial with respect to $s$. Consequently, $\deg Q_n = n,$ all its zeros are simple and lie in the open convex hull  of $\supp s$. Therefore, $\widehat{s}$ is a perfect system of one function.

In an attempt to construct general systems with abundant normal indices, E.M. Nikishin introduced in \cite{Nik} the concept of MT-system. Such systems are now named after him.

Let $\Delta_{\alpha}, \Delta_{\beta}$ be two intervals contained in the real line which do not intersect, or have at most a common end point, and $\sigma_{\alpha} \in {\mathcal{M}}(\Delta_{\alpha}), \sigma_{\beta} \in {\mathcal{M}}(\Delta_{\beta})$. Assume that $\widehat{\sigma}_{\beta} \in L_1(\sigma_{\alpha})$.  With these two measures we define a third one as follows (using the differential notation)
\[ d \la \sigma_{\alpha},\sigma_{\beta} \ra (x) = \widehat{\sigma}_{\beta}(x) d\sigma_{\alpha}(x) \in {\mathcal{M}(\Delta_{\alpha})}.
\]
Above, $\widehat{\sigma}_{\beta}$ denotes the Cauchy transform of the measure $\sigma_{\beta}$.  The more appropriate notation $\widehat{\sigma_{\beta}}$ causes space consumption and aesthetic inconveniences. We need to take consecutive products of measures; for example,
\[\la \sigma_{\gamma},  \sigma_{\alpha},\sigma_{\beta} \ra :=\la \sigma_{\gamma}, \la \sigma_{\alpha},\sigma_{\beta} \ra \ra. \]
 Here, we assume not only that $\widehat{\sigma}_{\beta} \in L_1(\sigma_{\alpha})$ but also $\la \sigma_{\alpha},\sigma_{\beta} \widehat{\ra} \in L_1(\sigma_{\gamma})$ where $\la \sigma_{\alpha},\sigma_{\beta} \widehat{\ra}$ denotes the Cauchy transform of $\la \sigma_{\alpha},\sigma_{\beta}  {\ra}$. Inductively, one defines products of a finite number of measures.

\begin{defi} \label{Nikishin} Take a collection  $\Delta_j, j=0,\ldots,m,$ of intervals such that, for each $j=0,\ldots,m-1$
\[ \Delta_j \cap \Delta_{j+1} = \emptyset, \qquad \mbox{or} \qquad \Delta_j \cap \Delta_{j+1} = \{x_{j,j+1}\},
\]
where $x_{j,j+1}$ is a single point. Let $(\sigma_0,\ldots,\sigma_m)$ be a system of measures such that $\mbox{Co}(\supp \sigma_j) = \Delta_j, \sigma_j \in {\mathcal{M}}(\Delta_j), j=0,\ldots,m,$  and
\begin{equation} \label{eq:autom}
\la \sigma_{j},\ldots,\sigma_k  {\ra} := \la \sigma_j,\la \sigma_{j+1},\ldots,\sigma_k\ra\ra\in {\mathcal{M}}(\Delta_j),  \qquad  0 \leq j < k\leq m.
\end{equation}
When $\Delta_j \cap \Delta_{j+1} = \{x_{j,j+1}\}$ we also assume that $x_{j,j+1}$ is not a mass point of either $\sigma_j$ or $\sigma_{j+1}$.
We say that $(s_0,\ldots,s_m) = {\mathcal{N}}(\sigma_0,\ldots,\sigma_m)$, where
\[ s_0 = \sigma_0, \quad s_1 = \la \sigma_0,\sigma_1 \ra, \ldots \quad , s_m = \la \sigma_0, \sigma_1,\ldots,\sigma_m  \ra
\]
is the Nikishin system of measures generated by $(\sigma_0,\ldots,\sigma_m)$.
\end{defi}

When we refer to a Nikishin system it is to be understood that all the assumptions made in Definition \ref{Nikishin} are satisfied. Condition (\ref{eq:autom}) is automatically fulfilled if $\Delta_j \cap \Delta_{j+1} = \emptyset$, provided $\la \sigma_{j+1},\ldots,\sigma_k  {\ra} \in {\mathcal{M}}(\Delta_{j+1})$. In particular, when $\Delta_j \cap \Delta_{j+1} = \emptyset, j=0,\ldots,m-1,$ (\ref{eq:autom}) is superfluous and we have the type of system originally defined by E. M. Nikishin.

Take $(s_0,\ldots,s_m) = {\mathcal{N}}(\sigma_0,\ldots,\sigma_m)$, where the moments of $s_0 = \sigma_0$ are finite.  Then, the moments of all $s_j, j=1,\ldots,m,$ are finite since $\lim_{x\to \infty, x \in \Delta_0} \la \sigma_1,\ldots,\sigma_j\widehat{\ra}(x) =0$ and $\la \sigma_1,\ldots,\sigma_j \widehat{\ra} \in L_1(\sigma_0)$.
Fix ${\bf n} \in {\mathbb{Z}}_+^{m+1}$ and consider the type II approximation of the Nikishin system of functions $(\widehat{s}_0,\ldots,\widehat{s}_m)$ relative to $\bf n$. It is easy to prove that
\[ \int x^{\nu} Q_{\bf n}(x) ds_j(x) =0, \qquad \nu = 0,\ldots,n_j -1,\qquad j=0,\ldots,m.
\]
If we denote
\[s_{j,k} = \la \sigma_j,\sigma_{j+1},\ldots,\sigma_k \ra,\qquad j < k, \qquad s_{j,j} = \la \sigma_j \ra = \sigma_j,
\]
the previous orthogonality relations may be rewritten as follows
\begin{equation}
\label{eq:a}
\int (p_0 (x) + \sum_{k=1}^m p_j(x)\widehat{s}_{1,k}(x))Q_{\bf n}(x) d \sigma_0(x) =0,
\end{equation}
where $p_0,\ldots,p_m$ are arbitrary polynomials such that $\deg p_k \leq n_k -1, k=0,\ldots,m.$

\begin{defi} \label{def:AT} A system of real continuous functions $u_0,\ldots,u_m$ defined on an interval $\Delta$ is called an AT-system on $\Delta$ for the multi-index ${\bf n} \in {\mathbb{Z}}_+^{m+1}$ if for any choice of real polynomials (that is, with real coefficients) $p_0,\ldots,p_m, \deg p_k \leq n_k -1,$ not all identically zero, the function
\[ \sum_{k=0}^m p_k(x) u_k(x)
\]
has at most $|{\bf n}| -1$ zeros on $\Delta$. If this is true for all ${\bf n} \in {\mathbb{Z}}_+^{m+1}$ we have an AT system on $\Delta$.
\end{defi}

In \cite{LF4} we proved that $(1,\widehat{s}_{1,1},\ldots,\widehat{s}_{1,m})$ forms an AT-system on any interval disjoint from $\Delta_1$ when the intervals $\Delta_j, j=1,\ldots,m,$ are bounded and $\Delta_j \cap \Delta_{j+1} = \emptyset, j=1,\ldots,m-1$. (For previous partial results, see also \cite{BBFG}, \cite{DrSt0},  \cite{LF1}, \cite{LIF}, \cite{GRS}, and \cite{Nik}.) Then, it easily follows from (\ref{eq:a}) that $\deg Q_{\bf n} = |{\bf n}|$, all its zeros are simple, and lie in the open convex hull  of $\supp \sigma_0$ (even if $\supp \sigma_0$ is unbounded or touches $\supp \sigma_1$). Thus, such Nikishin systems are type II perfect.  Moreover, we proved perfectness for mixed type Nikishin systems, having type I and type II as particular cases.

\subsection{Mixed type Nikishin systems.} Consider an $(m_2+1)\times (m_1+1)$ dimensional matrix of formal power expansions at $\infty$ ($m_2 +1$ rows and $m_1+1$ columns)
\[ {\mathbb{F}} = (f_{j.k}), \qquad f_{j,k}(z) = \sum_{n =0}^{\infty} \frac{f_{j,k,n}}{z^n},\qquad j=0,\ldots,m_2, \qquad k=0,\ldots,m_1.
\]
Fix a multi-index ${\bf n} = ({\bf n}_1;{\bf n}_2) \in {\mathbb{Z}}_+^{m_1 +1}\times {\mathbb{Z}}_+^{m_2 +1}$, such that $|{\bf n}_1| = |{\bf n}_2| +1$. We denote ${\bf n}_i =(n_{i,0},\ldots,n_{i,m_i}), i=1,2.$ There exists a  vector polynomial ${\mathbb{A}}_{\bf n}=(a_{{\bf n},0}, \ldots, a_{{\bf n},m_1})$,  such that
\begin{itemize}
\item[a)] ${\mathbb{A}}_{\bf n} \not \equiv {\bf 0}, \deg a_{{\bf n},k} \leq n_{1,k} -1, k=0,\ldots,m_1, $
\item[b)] $({\mathbb{F}}{\mathbb{A}}_{\bf n}^t - {\mathbb{B}}_{\bf n}^t)(z) = (\frac{A_{{\bf n},0}}{z^{n_{2,0}  +1}}+ \cdots, \ldots  ,\frac{A_{{\bf n},m_2}}{z^{n_{2,m_2}  +1}}+ \cdots)^t.$
\end{itemize}
for some $m_2 +1$ dimensional vector polynomial ${\mathbb{B}}_{\bf n}$ (the super-index $t$ means taking transpose and {\bf 0} denotes the zero vector). Finding ${\mathbb{A}}_{\bf n}$ reduces to solving a linear homogeneous system of $|{\bf n}_2|$ equations determined by the conditions b) on $|{\bf n}_1|$ unknowns (the total number of coefficients of the polynomials $a_{{\bf n},k}, k=0,\ldots,m_1$). Since $|{\bf n}_2| +1 = |{\bf n}_1|$ a non trivial solution exists.

\begin{defi}
\label{defiSor}
 A non zero vector ${\mathbb{A}}_{\bf n}$ satisfying a)-b) is called mixed type vector polynomial relative to ${\mathbb{F}}$ and ${\bf n} = ({\bf n}_1;{\bf n}_2)\in {\mathbb{Z}}_+^{m_1 +1}\times {\mathbb{Z}}_+^{m_2 +1}, |{\bf n}_1| = |{\bf n}_2| +1$. If $\deg a_{{\bf n},k} = n_{1,k} -1, k=0,\ldots,m_1$, the multi-index ${\bf n}$ is called mixed type normal. ${\mathbb{F}}$ is mixed type perfect when all multi-indices in ${\mathbb{Z}}_+^{m_1 +1}\times {\mathbb{Z}}_+^{m_2 +1}$ such that $|{\bf n}_1| = |{\bf n}_2| +1$ are normal.
\end{defi}

Mixed type systems were first introduced in \cite{Sor}.
This construction has as particular cases type I $(m_2 = 0)$ and type II $(m_1 =0)$ polynomials.
\vspace{0,2cm}

Let $S^{1} = (s^1_{0,0},\ldots,s^1_{0,m_1}) = {\mathcal{N}}(\sigma_0^1,\ldots,\sigma_{m_1}^1), S^{2} = (s^2_{0,0},\ldots,s^2_{0,m_2})= {\mathcal{N}}(\sigma_0^2,\ldots,\sigma_{m_2}^2), \sigma_0^1 = \sigma_0^2,$ be two given Nikishin systems generated by $m_1+1$ and $m_2+1$ measures, respectively. We underline  that both Nikishin systems stem from the same root measure $\sigma_0^1 = \sigma_0^2$.
Let us  introduce the row vectors
\[ {\mathbb{U}} =  (1,\widehat{s}^{2}_{1,1},\ldots,\widehat{s}^{2}_{1,m_2}),
\qquad {\mathbb{V }} =
(1,\widehat{s}^{1}_{1,1},\ldots,\widehat{s}^{1}_{1,m_1})
\]
and the $(m_2+1) \times (m_1 +1)$ dimensional  matrix function
\[ {\mathbb{W}} = {\mathbb{U}}^t {\mathbb{V}}.
\]
We say that the $S^1,S^2$ are compatible if asides from $\sigma_0^1 = \sigma_0^2$ we also have that
\[ {\mathbb{W}} \in L_1(\sigma_2^0);
\]
that is, each entry in ${\mathbb{W}}$ is integrable with respect to $\sigma_2^0$. In this case, we define the Markov (Stieltjes) type matrix function
\[ \widehat{\mathbb{S}}(z) = \int \frac{{\mathbb{W}}(x)d\sigma_0^2 (x)}{z-x}
\]
understanding that integration is carried out entry by entry on
the matrix $\mathbb{W}$. We say that $\widehat{\mathbb{S}}$ is a mixed type Nikishin system of functions.
We will study mixed type Nikishin systems.

For type II $(m_1 =0)$ we reduce the notation. Then, ${\bf n} = (n_0,\ldots,n_m)$, the vector function will be ${\bf f} = (\widehat{s}_0, \ldots, \widehat{s}_m)$, where $m = m_2$, and $(s_0,\ldots,s_m) = {\mathcal{N}}(\sigma_0,\ldots,\sigma_m)$. The mixed type polynomials ${\mathbb{A}}_{\bf n}$ will then be denoted $Q_{\bf n}$.

\subsection{Statement of the main results.}  Passing to the case of touching or unbounded supports in the generating measures of the Nikishin systems creates substantial technical difficulties.  The importance of the extended definition lies in the interesting examples it allows to consider. For example, take  ${\mathcal{N}}(\sigma_0,\sigma_1)$, where
\[ d\sigma_0(x) = e^{-x^{\lambda_1}}dx, \quad x \in [0,+\infty), \qquad d\sigma_1(x) = e^{x^{\lambda_2}}dx, \quad x \in (-\infty,0], \qquad \lambda_1,\lambda_2 > 0,
\]
or
\[ d\sigma_0(x) =  \frac{dx}{\sqrt{x(1-x)}}, \quad x \in [0,1], \qquad d\sigma_1(x) =  {dx} , \quad x \in [-1,0].
\]
These examples with classical weights, and their generalizations, have received considerable attention in brownian motion and random matrix theories (see, for example, \cite{BGS}, \cite{B}, \cite{DK}, \cite{DKV}, \cite{D}, and \cite{DF}) because of their use in describing important models. We would also like to mention here the recent paper $\cite{Sor2}$ dedicated to the study of the logarithmic asymptotic of type II multiple orthogonal polynomials for a Nikishin system of two unbounded measures which generalize Pollaczek polynomials.

In the results below, we have tried to keep the assumptions to a minimum.  For example, the proofs can be substantially simplified if one requires that the generating measures satisfy that whenever $\Delta_{j} \cap \Delta_{j+1} = \{x_{j,j+1}\}$ then $\widehat{\sigma}_j(x_{j,j+1}) \in{\mathbb{R}}$ and  $\widehat{\sigma}_{j+1}(x_{j,j+1}) \in{\mathbb{R}}$. Unfortunately, that condition is violated in the previous examples.

\begin{teo}
\label{teo:1}
Let $(s_{1,1},\ldots,s_{1,m}) = {\mathcal{N}}(\sigma_{1},\ldots,\sigma_m)$ be such that the first $2^k +1$  moments of $\sigma_k, k=1,\ldots,m,$ are finite.  Then, the system  $(1,\widehat{s}_{1,1},\ldots,\widehat{s}_{1,m})$ forms an AT-system on any interval $\Delta$ disjoint from $\Delta_1 = \mbox{\rm Co}(\supp \sigma_1)$. Moreover, for each ${\bf n} \in {\mathbb{Z}}_+^{m+1},$ and arbitrary polynomials with real coefficients $p_k, \deg p_k \leq n_k-1, k=0,\ldots,m,$ the linear form $p_0 + \sum_{k=1}^m p_k \widehat{s}_{1,k}, $ has at most $|{\bf n}|-1$ zeros in $\overline{\mathbb{C}} \setminus \Delta_1$, if it is not identically equal to zero.
\end{teo}

Definition \ref{defiSor} can be expressed in terms of orthogonality relations.  For that purpose, set
\[ {\mathcal{A}}_{{\bf n}}:= a_{{\bf n},0} + \sum_{k=1}^{m_1} a_{{\bf n},k} \widehat{s}_{1,k}^1.
\]

\begin{teo} \label{lem:orto} Let $S^1,S^2$ be two compatible Nikishin systems, such that the first $2^k +1$  moments of $\sigma_k^i, k=1,\ldots,m_i, i=1,2,$ are finite and all the moments of $\sigma_0^2 (= \sigma_0^1)$ are finite. Given ${\bf n} \in {\mathbb{Z}}_+^{m_1 +1}\times {\mathbb{Z}}_+^{m_2 +1}$ such that $|{\bf n}_1| = |{\bf n}_2| +1,$ we have
\begin{equation} \label{orto}
  \int {\mathcal{L}}_{{\bf n}_2}(x) {\mathcal{A}}_{{\bf n}}(x) d\sigma_0^2(x) = 0,
\end{equation}
for any linear form
\[ {\mathcal{L}}_{{\bf n}_2}(x) = p_0(x) + \sum_{j=1}^{m_2} p_j(x) \widehat{s}^2_{1,j}(x),
\]
where the $p_j, j=0,\ldots,m_2,$ denote arbitrary polynomials such that $\deg p_j \leq n_{2,j} -1$.    ${\mathcal{A}}_{{\bf n}}$ has exactly  $|{\bf n}_2|$ zeros in  ${\mathbb{C}} \setminus \mbox{\rm Co}(\supp \sigma_1^1)$, they are simple, and lie in the interior of $\mbox{\rm Co}(\supp \sigma_0^1)$.
\end{teo}

Because of (\ref{orto}), ${\mathbb{A}}_{\bf n}$ or ${\mathcal{A}}_{\bf n}$ are called mixed type orthogonal polynomials. An immediate consequence of Theorems \ref{teo:1} and \ref{lem:orto} is

\begin{teo}\label{teo:2} Let $S^1,S^2$ be two compatible Nikishin systems, such that the first $2^k +1$  moments of $\sigma_k^i, k=1,\ldots,m_i, i=1,2,$ are  finite  and all the moments of $\sigma_0^2 (= \sigma_0^1)$ are finite.  The matrix $\widehat{\mathbb{S}}$ is mixed type perfect. For each ${\bf n} \in {\mathbb{Z}}_+^{m_1 +1}\times {\mathbb{Z}}_+^{m_2 +1}, |{\bf n}_1| = |{\bf n}_2| +1,$ the vector polynomial ${\mathbb{A}}_{\bf n}$ is uniquely determined up to a constant factor.
\end{teo}

We say that a sequence of multi-indices $\Lambda \subset {\mathbb{Z}}^{m+1}, m \geq 0,$ is complete when it is totally ordered (componentwise) and $|\cdot|$ establishes a bijection between $\Lambda$ and $\mathbb{N}$ the set of natural numbers. From the previous results one has

\begin{cor} \label{biort} Let $S^1,S^2$ be two compatible Nikishin systems, such that the first $2^k +1$  moments of $\sigma_k^i, k=1,\ldots,m_i, i=1,2,$ are  finite  and all the moments of $\sigma_0^2 (= \sigma_0^1)$ are finite. Let $\Lambda_1 \subset {\mathbb{Z}}^{m_1+1}$ and $\Lambda_2 \subset {\mathbb{Z}}^{m_2+1}$ be two complete sequences of multi-indices. Then, there exist two sequences of linear forms $\{{\mathcal{Q}}_{{\bf n}_1}\}_{{\bf n}_1 \in \Lambda_1}$ and
$\{{\mathcal{P}}_{{\bf n}_2}\}_{{\bf n}_2 \in \Lambda_2}$ such that
\[ {\mathcal{Q}}_{{\bf n}_1} = q_{{\bf n}_1,0} + \sum_{k=1}^{m_1} q_{{\bf n}_1,k} \widehat{s}_{1,k}^1, \qquad \deg q_{{\bf n}_1,k} = n_{1,k} -1,\qquad  k=0,\ldots,m_1,
\]
\[ {\mathcal{P}}_{{\bf n}_2} = p_{{\bf n}_2,0} + \sum_{k=1}^{m_2} p_{{\bf n}_2,k} \widehat{s}_{1,k}^2, \qquad \deg p_{{\bf n}_2,k} = n_{2,k} -1,\qquad  k=0,\ldots,m_2,
\]
each  ${\mathcal{Q}}_{{\bf n}_1}$ and ${\mathcal{P}}_{{\bf n}_2}$ is uniquely determined except for a constant factor,
\begin{equation} \label{bi1} \int {\mathcal{Q}}_{{\bf n}_1}(x) {\mathcal{P}}_{{\bf n}_2}(x) d \sigma_0^2(x) =0,\qquad |{\bf n}_1| \neq |{\bf n}_2| ,
\end{equation}
and
\begin{equation} \label{bi2} \int {\mathcal{Q}}_{{\bf n}_1}(x) {\mathcal{P}}_{{\bf n}_2}(x) d \sigma_0^2(x) \neq 0,\qquad |{\bf n}_1| = |{\bf n}_2|.
\end{equation}
${\mathcal{Q}}_{{\bf n}_1}$ has exactly $|{\bf n}_1|-1$ zeros in ${\mathbb{C}} \setminus \mbox{\rm Co}(\supp(\sigma_1^1))$, they are all simple, and lie in the interior of $\mbox{\rm Co}(\supp(\sigma_0^1))$. Likewise, ${\mathcal{P}}_{{\bf n}_2}$ has exactly $|{\bf n}_2|-1$ zeros in ${\mathbb{C}} \setminus \mbox{\rm Co}(\supp(\sigma_1^2))$, they are all simple, and lie in the interior of $\mbox{\rm Co}(\supp(\sigma_0^2))$.
\end{cor}

A very general treatment on the construction and algebraic properties of general bi-orthogonal sequences of generalized polynomials such as $(\{{\mathcal{Q}}_{{\bf n}_1}\}_{{\bf n}_1 \in \Lambda_1}, \{{\mathcal{P}}_{{\bf n}_2}\}_{{\bf n}_2 \in \Lambda_2})$ may be found in \cite{AFM}. In particular, their recursion and Christoffel-Darboux type formulas are derived as well as their connection to the multi-component 2D Toda lattice hierarchy is studied. Christoffel-Darboux type formulas are essential in the description of the correlation function in brownian motion and random matrix models (see \cite{BGS}, \cite{B}, \cite{DK}, \cite{DKV}, \cite{D}, and \cite{DF}).

If we restrict our attention to type II approximation, the location of the zeros of $Q_{\bf n} (=  {\mathcal A}_{\bf n})$ on the support of $\sigma_0$ together with \cite[Theorem 1]{Bus} allow us to obtain an analogue of the Stieltjes theorem (see \cite{Sti}) on the convergence of diagonal Pad\'e approximation.

\begin{cor}
\label{cor:1} Let $(s_{0},\ldots,s_{m}) = {\mathcal{N}}(\sigma_0,\ldots,\sigma_m)$ be such that all the moments of $\sigma_0$ are finite, and for each $k=1,\ldots,m,$ the first $2^k +1$ moments of $\sigma_k$ are finite. Let $\Lambda \subset {\mathbb{Z}}_+^{m+1}$ be a sequence of distinct multi-indices such that
\[ n_j \geq \frac{|{\bf n}|}{m+1}  - c  ,  \qquad j=0,\ldots,m,
\]
for some   constant $c > 0$.  Assume that one of the following conditions is verified:
\begin{itemize}
\item[a)] $\Delta_0 \cap \Delta_1 = \emptyset$ and either $\Delta_0$ or $\Delta_1$ is bounded.
\item[b)] $\Delta_0$ is unbounded and
\[ \sum_{n=0}^{\infty} \frac{1}{|c_n|^{1/2n}} = \infty, \qquad c_n = \int x^n d\sigma_0(x).
\]
\item[c)] $\Delta_0 \cap \Delta_1 = \{x_{0,1}\}$ and
\[ \sum_{n=0}^{\infty} \frac{1}{|d_n|^{1/2n}} = \infty, \qquad d_n = \int \frac{ d\sigma_0(x)}{(x-x_{0,1})^n}.
\]
\end{itemize}
Then,
\[ \lim_{{\bf n} \in \Lambda} \frac{P_{{\bf n},k}}{Q_{\bf n}} = \widehat{s}_{0,k} , \qquad k=0,\ldots,m
\]
uniformly on each compact ${\mathcal{K}} \subset \overline{\mathbb{C}} \setminus \Delta_0$. When $a)$ occurs, the convergence has geometric rate.
\end{cor}

For measures with unbounded support, this result extends \cite[Corollary 1]{Bus} to a more general class of multi-indices allowing also touching supports. When the supports are bounded and non intersecting, \cite[Corollary 1.1]{LF4} gives more information.

\section{Auxiliary results}

Generally speaking, the underlying arguments in the proofs of the theorems stated above are the same as for similar ones contained in \cite{LF4}. Justifying that they are feasible in the present setting is technically more difficult and that is our task in this section. Whenever needed, we will rely on \cite{LF4} to avoid unnecessary repetitions.

The following result, similar to \cite[Lemma 2.2]{LF4}, covers any system satisfying Definition \ref{Nikishin}.

\bl
\label{lem:3}
Let $(s_{1,1},\ldots,s_{1,m}) =
{\mathcal{N}}(\sigma_1,\ldots,\sigma_m)$ and ${\bf n} \in
{\mathbb{Z}}_+^{m+1}$ be given. Consider the linear form
\[ {\mathcal{L}}_{\bf n}  = p_0 + \sum_{k=1}^m p_k \widehat{s}_{1,k}, \quad \deg p_k \leq n_k-1, \quad k=0,\ldots,m,
\]
where the polynomials $p_k$ have real coefficients. Assume that $n_0 = \max\{n_0,n_1-1,\ldots,n_m-1\}$. If
${\mathcal{L}}_{\bf n}$ had at least $|{\bf n}|$ zeros in $ {\mathbb{C}} \setminus \Delta_1$ the reduced form
$p_1 + \sum_{k=2}^m p_k \widehat{s}_{2,k}$ would have at least $|{\bf n}| - n_0$ zeros in $\Delta_1 \subset {\mathbb{C}} \setminus \Delta_2$.
\el

{\bf Proof.} The function ${\mathcal{L}}_{\bf n}$ is symmetric with respect to the real line, ${\mathcal{L}}_{\bf n}(\overline{z}) = \overline{{\mathcal{L}}_{\bf n}(z)}$; therefore, its zeros come in conjugate pairs. Thus, if ${\mathcal{L}}_{\bf n}$ has at least $|{\bf n}|$ zeros in $ {\mathbb{C}} \setminus \Delta_1$, there exists a polynomial $w_{\bf n}, \deg w_{\bf n}
\geq |{\bf n}|,$  with real coefficients and zeros contained in
${\mathbb{C}}\setminus \Delta_1$ such that ${\mathcal{L}}_{\bf
n}/w_{\bf n} \in {\mathcal{H}}({\mathbb{C}}\setminus \Delta_1)$. This function has a zero of order $ \geq |{\bf n}| - n_0 +1$
at $\infty$. Consequently,
\[ \frac{z^{\nu} {\mathcal{L}}_{\bf n}}{w_{\bf n}}  \in {\mathcal{H}}( {\mathbb{C}} \setminus
\Delta_1),
\]
and, for all $ \nu=0,\ldots,|{\bf n}| - n_0 -1,$
\[ \deg z^{\nu}p_k \leq |{\bf n}| -n_0 - 1 + n_k -1 \left\{
\begin{array}{cc}
= |{\bf n}| -2, & k=0 \\
\leq |{\bf n}| -1, & k=1,\ldots,m.
\end{array}
\right.
\]

Let $\Gamma$ be a closed simple curve which surrounds all the zeros of $w_{\bf n}$ and leaves $\Delta_1$ outside. Since $\frac{z^{\nu} {\mathcal{L}}_{\bf n}}{w_{\bf n}}$ is analytic inside $\Gamma$ and the functions $\frac{z^{\nu}p_k}{w_{\bf
n}}$ are analytic outside, with a zero of second order at infinity if $k=0$ and order $1$ when $k=1,\ldots,m$, integrating along $\Gamma$ we have ($\widehat{s}_{2,1} \equiv 1$)
\[ 0 = \int_{\Gamma} \frac{z^{\nu} {\mathcal{L}}_{\bf n}(z)dz}{w_{\bf n}(z)} =  \sum_{k=1}^m \int_{\Gamma} \frac{z^{\nu}  p_k (z)}{w_{\bf n}(z)}
\int \frac{d{s}_{1,k}(x)}{z-x} dz = \sum_{k=1}^m \int \widehat{s}_{2,k}(x) \int_{\Gamma} \frac{z^{\nu}  p_k (z)}{w_{\bf n}(z)(z-x)} dz d \sigma_{1}(x) =
\]
\[ -2\pi i \int x^{\nu}(p_1(x) + \sum_{k=2}^m p_k(x) \widehat{s}_{2,k}(x)) \frac{d\sigma_1(x)}{w_{\bf n}(x)}, \qquad \nu= 0,\ldots,|{\bf n}| -n_0 -1.
\]

These orthogonality relations imply that $p_1 + \sum_{k=2}^m p_k
\widehat{s}_{2,k}$ has at least $|{\bf n}| - n_0$ sign changes in the interior of $\Delta_1$.  In fact, if there were at most $|{\bf n}| - n_0 -1$ sign changes one can easily construct a polynomial $p$ of degree $\leq |{\bf n}| - n_0 -1$ such that $p(p_1 + \sum_{k=2}^m p_k
\widehat{s}_{2,k})$ does not change sign on $\Delta_1$ which contradicts the orthogonality relations. Therefore, in the interior of $\Delta_1 \subset {\mathbb{C}} \setminus \Delta_2$, the reduced form would have the number of zeros claimed.  \hfill $\Box$ \vspace{0,2cm}

\bl
\label{lem:3a}
Let $(s_{1,1},\ldots,s_{1,m}) =
{\mathcal{N}}(\sigma_1,\ldots,\sigma_m)$ be given. Then $(1,\widehat{s}_{1,1},\ldots,\widehat{s}_{1,m})$ forms an AT system with respect to any multi-index ${\bf n} \in
{\mathbb{Z}}_+^{m+1}(\bullet)$, where
\[  {\mathbb{Z}}_+^{m+1}(\bullet) = \{{\bf n} \in {\mathbb{Z}}_+^{m+1}:  n_0 \geq \cdots \geq n_m\},
\]
on any interval $\Delta$ disjoint from $\Delta_1 = \mbox{\rm Co}(\supp \sigma_1)$. Moreover, the linear form
\[ {\mathcal{L}}_{\bf n} = p_0 + \sum_{k=1}^m p_k \widehat{s}_{1,k}, \quad \deg p_k \leq n_k-1, \quad k=0,\ldots,m,
\]
has at most $|{\bf n}| -1$ zeros in ${\mathbb{C}} \setminus \Delta_1$ if it is not identically equal to zero.
\el

{\bf Proof.} Take ${\bf n} \in {\mathbb{Z}}_+^{m+1}(\bullet)$ and some ${\mathcal{L}}_{\bf n}$ not identically zero. Suppose that $p_j$ is the polynomial in the form ${\mathcal{L}}_{\bf n}$ with largest subindex  which is not identically zero. Assume that ${\mathcal{L}}_{\bf n}$ has at least $|{\bf n}|$ zeros in ${\mathbb{C}} \setminus \Delta_1$. Applying $j$ times Lemma \ref{lem:3}, it follows that $p_j$ has at least $|{\bf n}| - (n_0+\cdots+n_{j-1}) \geq n_j$ zeros in $\Delta_{j}$. This is absurd because $\deg p_j \leq n_j -1 $.   \hfill $\Box$ \vspace{0,2cm}

Lemma \ref{lem:3a} settles Theorem \ref{teo:1} if we restrict our attention to multi-indices with decreasing components. For multi-indices with decreasing components we have only needed the existence of the first moment; that is, that the measures be finite. The proof of the general result reduces to that of multi-indices with decreasing components after carrying out intricate transformations of the generating system of measures which do require the existence of some moments of higher order.

An important ingredient in these transformations consists in an appropriate representation for  the reciprocal of the Cauchy transform of measures. For measures with bounded support, this is known and was used in \cite{LF4}. In the unbounded case, we could not find a proper reference of the results and properties needed, though they are known to the specialists. For completeness, we include a brief account.

Let $s$ be a measure with constant sign, supported on  an interval $\Delta \subset {\mathbb{R}}$, and finite moments
\[c_n = \int x^n d s(x), \qquad n \in {\mathbb{Z}}_+,
\]
It is well known that
\[ \widehat{s}(z) \sim \sum_{n=0}^{\infty} \frac{c_n}{z^{n+1}} .
\]
If $\supp s$ is bounded,  we can write $=$ instead of $\sim$ for all $z$ in a neighborhood of $\infty$.
By $\sim$ we mean that for each $n \in {\mathbb{Z}}_+$
\begin{equation} \label{eq:3} \lim_{z \to \infty} z^{n+1}(\widehat{s}(z) - \sum_{k=0}^{n-1} \frac{c_k}{z^{k+1}}) = c_{n}.
\end{equation}
The limit is taken along any sector which does not intersect $\Delta$ and limiting rays not parallel to $\Delta$.

When the moments $c_0,\ldots,c_{n-1}$ of $s$ are finite,  we have
\begin{equation} \label{resto} \int \frac{d s(x)}{z-x} = \frac{c_0}{z} + \cdots + \frac{c_{n-1}}{z^{n}} + \int \frac{x^{n} ds(x)}{z^{n}(z-x)},
\end{equation}
which is valid in the complement  of $\supp s$. Consequently,
\[ z^{n+1}(\widehat{s}(z) - \sum_{k=0}^{n-1} \frac{c_k}{z^{k+1}}) = \int  \frac{z}{ z-x} x^{n} ds(x)
\]
($\sum_{k=0}^{-1}$denotes an empty sum). Due to Lebesgue's dominated  convergence theorem it follows that $c_n$  is   finite if and only if the limit in (\ref{eq:3}) is finite, coinciding  the value of the  limit and the moment.

\bl \label{momentos} Let $s \in {\mathcal{M}}(\Delta)$, where $\Delta$ is a half line of ${\mathbb{R}}$, and its first two moments $ c_0, c_1$ are finite. Then, there exists  a measure $\tau$ with sign opposite to that of $s$ and support contained in $\Delta$ such that $\int (1 + |x|)^{-1} d|\tau|(x)  < \infty$ and
\begin{equation}\label{eq:1} \frac{1}{\widehat{s}(z)}= \frac{1}{c_0}z - \frac{c_1}{c_0^2} + \widehat{\tau}(z).
\end{equation}
If $a$ denotes the finite end point of $\Delta$, then $|\widehat{\tau}(a)| < \infty$; in particular, $a$ is not a mass point of $\tau$.
When the first $n+3, n \geq 0,$ moments of $s$ are finite then the first $n+1$  moments of $\tau$ are finite. Thus, $\tau$ is finite if the first three moments of $s$ are finite.
\el

{\bf Proof.} Without loss of generality, we can assume that $s$ is a positive measure and $\Delta = [0,+\infty)$. Notice that $1/z\widehat{s}(z)$ is analytic in ${\mathbb{C}} \setminus [0,\infty)$.
Writing $z = t +iy, t,y \in {\mathbb{R}}$, we have
\[ \Im \frac{1}{z\widehat{s}(z)} = \Im \frac{\overline{z} \widehat{s}(\overline{z})}{|z\widehat{s}(z)|^2} =
\frac{y}{|z\widehat{s}(z)|^2}\int \frac{xds(x)}{|z-x|^2} > 0 \quad \mbox{for} \quad y >0, \quad \mbox{and}\quad  \frac{1}{t\widehat{s}(t)} > 0 \quad \mbox{for} \quad t <0,
\]
where $\Im$ denotes the imaginary part. Thus, by \cite[Theorem A.4]{KN}, there exists a positive measure $\tau_1$ supported in $[0,+\infty)$ with $\int (1+x)^{-1}{d\tau_1 (x)} < \infty$ and a non negative constant $d_{-2}$ such that
\[ \frac{1}{z\widehat{s}(z)} = d_{-2} + \int \frac{d\tau_1 (x)}{x-z} = d_{-2} + \int \frac{1+x}{x-z}\frac{d\tau_1 (x)}{1+x}.
\]
Since $\lim_{t \to -\infty} t\widehat{s}(t) = c_0$, using Lebesgue's dominated convergence theorem, we have $d_{-2} = 1/c_0$.

From (\ref{resto})  with $n=0$, we obtain
\[  \frac{1}{\widehat{s}(z)} - d_{-2} z = \frac{z(1- d_{-2}z \widehat{s}(z))}{z \widehat{s}(z)} = - \frac{d_{-2}  \int\frac{z }{z-x} xds(x)}{z \widehat{s}(z)} = - \int \frac{z d\tau_1 (x)}{z-x}  ,
\]
and, by Lebesgue's monotone convergence theorem, it follows  that
\[ \lim_{t\to -\infty} \frac{1}{\widehat{s}(t)} - d_{-2} t = -\frac{c_1}{c_0^2} := d_{-1} = - \int d\tau_1 (x)\in {\mathbb{R}}.
\]
Consequently,
\[ \frac{1}{\widehat{s}(z)} = \frac{z}{c_0} - \frac{c_1}{c_0^2} + \int \frac{d \tau(x)}{z -x},
\]
where $d\tau (x) = - x d\tau_1 (x)$. Therefore (\ref{eq:1}) holds since
\[ \int \frac{d|\tau| (x)}{1+x} < \int d\tau_1 (x)  < \infty.
\]

Notice that the function $\widehat{s}(t)$ is negative and decreasing for $t<0$. Therefore,
$\lim_{t \to 0-} 1/\widehat{s}(t)$ exists and is either $0$ or some finite negative number.
Consequently, $\lim_{t \to 0-} \widehat{\tau}(t)$ exists and it is a finite number. Using Lebesgue's monotone convergence theorem it follows that $1/x \in L_1(\tau)$, or what is the same, $|\widehat{\tau}(0)| < \infty$ as we needed to prove.

Now, let us assume that $c_0,\ldots,c_{n+2}$ are finite. Define consecutively numbers $d_n \in {\mathbb{R}}, n \in \mathbb{Z}_+$, through the following  triangular  system of equations (later we prove that they are moments of $\tau$)
\begin{equation} \label{eq:4}
\begin{array}{ccl}
1 & = & d_{-2} c_0 \\
0 & = & d_{-2} c_1 + d_{-1} c_0\\
0 & = & d_{-2} c_2 + d_{-1} c_1 + d_0 c_0\\
\vdots & = & \vdots  \\
0 & = & d_{-2} c_{n+2} + d_{-1} c_{n+1} + \cdots + d_n c_0\,\,.
\end{array}
\end{equation}
In particular, we have
\[ d_{-2} = \frac{1}{c_0}, \qquad d_{-1} = \frac{-c_1}{c_0^2}, \qquad d_0 = \frac{-c_2}{c_0^2} +\frac{c_1^2}{c_0^3}, \]
and so on. The first two values coincide with the definition given above for $d_{-2},d_{-1}$.

We will show that
\[ \lim_{z \to \infty} z^{k+1}(\widehat{\tau}(z) - \sum_{j=0}^{k-1} \frac{d_j}{z^{j+1}}) = d_{k}, \qquad k=0,\ldots,n.
\]
As noted above, this implies that $d_k = \int x^k d \tau(x), k=0,\ldots,n$.

Let us prove that
\[ \lim_{z\to \infty} z \widehat{\tau}(z) = \lim_{z\to \infty} z (\frac{1}{\widehat{s}(z)} - (d_{-2} z + d_{-1})) = d_0. \]
In fact, using (\ref{resto}), the first three equations in (\ref{eq:4}),  and Lebesgue's dominated convergence theorem, it follows that
\[ z (\frac{1}{\widehat{s}(z)} - (d_{-2} z + d_{-1})) = z^2 \frac{1 -  (d_{-2} z + d_{-1})\widehat{s}(z) }{z\widehat{s}(z)} =
\]
\[ z^2  \frac{1 -  (d_{-2} z + d_{-1})(\frac{c_0}{z} + \frac{c_1}{z^2} + \frac{c_2}{z^3} + \int \frac{x^3 d s(x)}{z^3(z-x)})}{z\widehat{s}(z)} =
\frac{d_0c_0 - \frac{d_{-1}c_2}{z} - \frac{d_{-2} z - d_{-1}}{z} \int \frac{x^3 ds(x)}{z-x}}{z\widehat{s}(z)} =
\]
\[ \frac{d_0c_0 + {\mathcal{O}}(1/z) + { \mathcal{O}}(1)\int \frac{x^3 ds(x)}{z-x} }{z\widehat{s}(z)} \longrightarrow_{z \to \infty}  d_0.
\]
In general, for $k=1,\ldots,n$, on account of (\ref{resto}), the first $k+3$ equation in (\ref{eq:4}), and Lebesgue's dominated convergence theorem, we obtain
\[ z^{k+1} (\frac{1}{\widehat{s}(z)} - (d_{-2} z + d_{-1} + \sum_{j=0}^{k-1} \frac{d_j}{z^{j+1}})) =
z^{k+2} \frac{1 -   (\sum_{j=-2}^{k-1} \frac{d_j}{z^{j+1}}) (\sum_{i=0}^{k+2} \frac{c_i}{z^{i+1}} + \int \frac{x^{k+3} ds(x)}{z^{k+3}(z-x)})} {z\widehat{s}(z)} =
\]
\[ \frac{d_k c_0 + {\mathcal{O}}(1/z) + { \mathcal{O}}(1)\int \frac{x^{k+3} ds(x)}{z-x}}{z\widehat{s}(z)} \longrightarrow_{z \to \infty}  d_k,
\]
as we needed to prove. \hfill $\Box$
\vspace{0,2cm}

In the sequel, we write
\[
{1}/{\widehat{s}(z)}={\ell}(z)+ \widehat{\tau}(z),
\]
where $\ell$ denotes the first degree polynomial in the decomposition (\ref{eq:1}). For
convenience, we call $\tau$  the inverse measure of $s.$ Such
measures will appear frequently in our reasonings, so we will fix a
notation to distinguish them. They will always refer to inverses of
measures denoted with $s$ and will carry over to them the
corresponding sub-indices. The same goes for the  polynomials
$\ell$. For instance, if $s_{\alpha,\beta} = \langle
\sigma_{\alpha},\sigma_{\beta} \rangle $, then
\[
{1}/{\widehat{s}_{\alpha,\beta}(z)}  ={\ell}_{\alpha,\beta}(z)+
\widehat{\tau}_{\alpha,\beta}(z).
\]
Sometimes we write $\langle
\sigma_{\alpha},\sigma_{\beta} \widehat{\rangle}$ in place of $\widehat{s}_{\alpha,\beta}$. This is specially useful later on where we need the Cauchy transforms of complicated expressions of products of measures for which we do not have a short hand notation. Since $s_{\alpha,\alpha} = \sigma_{\alpha}$, we also write
\[
{1}/{\widehat{\sigma}_{\alpha}(z)} ={\ell}_{\alpha,\alpha}(z)+
\widehat{\tau}_{\alpha,\alpha}(z).
\]

\begin{rem} \label{rem:1} Notice that from Lemma \ref{momentos} the inversion of $\tau$ requires the existence of its first two moments, if we wish to obtain for $1/\widehat{\tau}$ a formula like (\ref{eq:1}). By the same lemma this is guaranteed if the first $4$ moments of $s$ are finite. In general, $k$ consecutive inversions of a measure $s$ requires that its first $2^k$ moments be finite. On the other hand, the inversion of $s_{\alpha,\beta}$ requires two moments of $s_{\alpha,\beta}$, but this is true if the first two moments of $\sigma_{\alpha}$ are finite since $\lim_{z \to \infty}\widehat{\sigma}_{\beta}(z) = 0$ and  $\widehat{\sigma}_{\beta} \in L_1(\sigma_{\alpha})$. Therefore, $s_{\alpha,\beta}$ may be inverted $k$ times if the first $2^k$ moments of $\sigma_{\alpha}$ exist (or $2^k +1$ moments if we want that the final measure be finite).
\end{rem}

To deduce the reduction formulas, we need to apply Cauchy's theorem, Cauchy's integral formula, and Fubini's theorem, integrating along curves which surround intervals $\Delta_j$ and pass between consecutive intervals $\Delta_j, \Delta_{j+1}$. When these intervals are unbounded and/or touch, the curve has to go through infinity and the intersection point producing singularities in the kernels of the integrals. We overcome these difficulties considering principal values of improper integrals.

Fix $\theta \in (0,\pi), a \in {\mathbb{R}},$ and $\varepsilon, R > 0$.  Set
\[\Gamma_{\theta} = \{z= a + te^{i\theta}, t \geq 0\} \cup \{z =a+te^{-i\theta}, t \geq 0\},
\]
and
\[ \Gamma_{\epsilon,R,\theta} = \{z = a + te^{i\theta}: t \in [\varepsilon,R]\} \cup \{z = a+te^{-i\theta}: t \in [\varepsilon,R]\}.
\]
We consider that $\Gamma_{\theta}$ is oriented so that $(a,+\infty)$ is to the right of $\Gamma_{\theta}$ as we walk along $\Gamma_{\theta}$. $\Gamma_{\epsilon,R,\theta}$ has the orientation induced by $\Gamma_{\theta}$. Let $f: \Gamma_{\theta} \setminus \{a\} \longrightarrow {\mathbb{C}}$. We denote
\[ \int_{\Gamma_{\theta}} f(z) dz = \lim_{\varepsilon \to 0, R \to \infty} \int_{\Gamma_{\epsilon,R,\theta}} f(z) dz
\]
whenever the limit on the right hand exists and is finite. Set
\[ C_{R,\theta,1} = \{z= a + Re^{it}: \theta \leq t \leq 2\pi - \theta\}, \qquad C_{R,\theta,2} = \{z= a+Re^{it}: -\theta \leq t \leq \theta\}.
\]
Analogously,
\[ C_{\varepsilon,\theta,1} = \{z= a + \varepsilon e^{it}: \theta \leq t \leq 2\pi - \theta\}, \qquad C_{\varepsilon,\theta,2} = \{z= a+ \varepsilon e^{it}: -\theta \leq t \leq \theta\}.
\]
We assume that $C_{\varepsilon,\theta,1} \cup \Gamma_{\varepsilon,R,\theta} \cup C_{R,\theta,1}$ is oriented positively and $C_{\varepsilon,\theta,2} \cup \Gamma_{\varepsilon,R,\theta} \cup C_{R,\theta,2}$ negatively.

The next result substitutes Lemma 2.1 of \cite{LF4}.

\bl \label{lem:2} Fix $\theta, 0 < \theta < \pi.$ Let $h \in {\mathcal{H}}({\mathbb{C}}\setminus (-\infty,a])$ be such that
\begin{itemize}
\item[a)] $\sup_{z \in C_{R,\theta,2}} |h(z)| = {\mathcal{O}}(1/R^2), R \to \infty,$
\item[b)] $\sup_{z \in C_{\varepsilon,\theta,2}} |(z-a)h(z)| = {o}(1), \varepsilon \to 0,$
\end{itemize}
then
\begin{equation} \label{eq:2.1}
0 = \int _{\Gamma_{\theta}} h(\zeta) d \zeta .
\end{equation}
Let $w \in {\mathcal{H}}({\mathbb{C}}\setminus [a,+\infty))$ be such that
\begin{itemize}
\item[c)] $\sup_{z \in C_{R,\theta,1}} |w(z)| = {\mathcal{O}}(1/R), R \to \infty,$
\item[d)] $\sup_{z \in C_{\varepsilon,\theta,1}} |(z-a)w(z)| = {o}(1), \varepsilon \to 0,$
\end{itemize}
then, for any compact set ${\mathcal{K}} \subset {\mathbb{C}} \setminus [a,+\infty)$ to the left of $\Gamma_{\theta}$, \begin{equation} \label{eq:2.2}
w(z) = \frac{1}{2\pi i} \int _{\Gamma_{\theta}} \frac{w(\zeta) d \zeta}{\zeta - z}, \qquad z \in {\mathcal{K}}.
\end{equation}
Assume that $s \in {\mathcal{M}}[a,+\infty)$ is such that $|\widehat{s}(a)| < \infty$, and $g \in {\mathcal{H}}({\mathbb{C}} \setminus (-\infty,a])$ with constant sign on $(a,+\infty)$ verifies
\begin{itemize}
\item[e)] $\sup_{z \in C_{R,\theta,2}} |g(z)| = {\mathcal{O}}(1), R \to \infty,$
\item[f)] $\sup_{z \in C_{\varepsilon,\theta,2}} |(z-a)g(z)| = {o}(1), \varepsilon \to 0,$
\end{itemize}
then, $g d s \in {\mathcal{M}}[a,+\infty)$ and
\begin{equation} \label{eq:2.3}
\frac{1}{2\pi i}\int_{\Gamma_{\theta}} \frac{g(\zeta) \widehat{s}(\zeta)}{\zeta -z} d \zeta = \int \frac{g(x) d s(x)}{z-x}
\end{equation}
if the integral on the left hand is finite.
\el

{\bf Proof.} Without loss of generality, we assume that $a =0$. By Cauchy's theorem
\[ 0 = \int_{C_{\varepsilon,\theta,2} \cup \Gamma_{\varepsilon,R,\theta} \cup C_{R,\theta,2}} h(\zeta) d\zeta.
\]
Since
\[ |\int_{C_{\varepsilon,\theta,2}} h(\zeta) d\zeta| \leq 2\pi \sup_{\zeta \in C_{\varepsilon,\theta,2}} |\zeta h(\zeta)| \quad \mbox{and} \quad |\int_{C_{R,\theta,2}} h(\zeta) d\zeta| \leq 2\pi R\sup_{\zeta \in C_{R,\theta,1}} |h(\zeta)|,
\]
formula (\ref{eq:2.1}) follows using a) and b).

Next, given a compact set ${\mathcal{K}} \subset {\mathbb{C}} \setminus [0,+\infty)$ to the left of $\Gamma_{\theta}$ for all $\varepsilon$  sufficiently small and $R$ sufficiently large, the closed curve $C_{\varepsilon,\theta,1} \cup \Gamma_{\varepsilon,R,\theta} \cup C_{R,\theta,1}$ has winding number 1 with respect to all points in ${\mathcal{K}}$.
From Cauchy's integral formula, we have that
\[ w(z) = \frac{1}{2\pi i} \int_{C_{\varepsilon,\theta,1} \cup \Gamma_{\varepsilon,R,\theta} \cup C_{R,\theta,1}} \frac{w(\zeta) d\zeta}{\zeta - z}.
\]
But
\[ |\frac{1}{2\pi i} \int_{C_{\varepsilon,\theta,1}} \frac{w(\zeta) d\zeta}{\zeta - z}| \leq \frac{1}{||z| - \varepsilon|} \sup_{\zeta \in C_{\varepsilon,\theta,1}} |\zeta w(\zeta)|
\]
and
\[ |\frac{1}{2\pi i} \int_{C_{R,\theta,1}} \frac{w(\zeta) d\zeta}{\zeta - z} | \leq \frac{R}{|R - |z||}\sup_{\zeta \in C_{R,\theta,1}} |w(\zeta)|.
\]
therefore, (\ref{eq:2.2}) follows from c) and d).

That $g d s \in {\mathcal{M}}[a,+\infty)$ follows easily from e), f), and $|\widehat{s}(0)| < \infty\,\, (1/x \in L_1(s))$. Let us prove (\ref{eq:2.3}). By definition
\[ \int_{\Gamma_{\theta}} \frac{g(\zeta) \widehat{s}(\zeta)}{\zeta -z} d \zeta = \lim_{\varepsilon \to 0, R \to \infty} \int_{\Gamma_{\varepsilon,R,\theta}} \frac{g(\zeta) \widehat{s}(\zeta)}{\zeta -z} d \zeta,
\]
and
\[\int_{\Gamma_{\varepsilon,R,\theta}} \frac{g(\zeta) \widehat{s}(\zeta)}{\zeta -z} d \zeta = \int_{\Gamma_{\varepsilon,R,\theta}}   \int_{[0,2\varepsilon)\cup[2\varepsilon,R/2]\cup(R/2,+\infty)}  \frac{d s(x) }{\zeta - x}\frac{g(\zeta)d \zeta}{\zeta - z}.
\]
Let us study the different integrals.

Notice that
\[\int_{\Gamma_{\varepsilon,R,\theta}}   \int_{[0,2\varepsilon)\cup(R/2,+\infty)}  \frac{d s(x) }{\zeta - x}\frac{g(\zeta)d \zeta}{\zeta - z}  =
\]
\[\int_{\Gamma_{\varepsilon,1,\theta}}   \int_{[0,2\varepsilon)\cup(R/2,+\infty)}  \frac{d s(x) }{\zeta - x}\frac{g(\zeta)d \zeta}{\zeta - z} + \int_{\Gamma_{1,R,\theta}}   \int_{[0,2\varepsilon)\cup(R/2,+\infty)}  \frac{d s(x) }{\zeta - x}\frac{g(\zeta)d \zeta}{\zeta - z}.
\]
We have
\[ |\int_{\Gamma_{1,R,\theta}}   \int_{[0,2\varepsilon)\cup(R/2,+\infty)}  \frac{d s(x) }{\zeta - x}\frac{g(\zeta)d \zeta}{\zeta - z}| \leq
\]
\[ \sup_{\zeta \in \Gamma_{1,R,\theta}} |g(\zeta)| \int_{\Gamma_{1,R,\theta}}   \int_{[0,2\varepsilon)\cup(R/2,+\infty)} \frac{|\zeta|}{|\zeta - x|} d| s|(x)  \frac{|\zeta|}{|\zeta - z|}\frac{|d \zeta|}{|\zeta|^2} \leq
\]
\[ \frac{M}{\sin \theta} |s|([0,2\varepsilon)\cup(R/2,+\infty)) \longrightarrow 0, \quad \varepsilon \to 0, \quad R \to \infty.
\]
Here, we used that $\sup_{\zeta \in \Gamma_{1,R,\theta}} |g(\zeta)|$ is uniformly bounded with respect to $R$, ${|\zeta|}/{|\zeta - x|} \leq 1/\sin \theta$, $\sup_{z \in {\mathcal{K}}, \zeta \in \Gamma_{\theta}}\frac{|\zeta|}{|\zeta - z|}$ is finite, $\int_1^{\infty} dt/t^2 =1$, and $\lim_{\varepsilon \to 0, R \to \infty} |s|([0,2\varepsilon)\cup(R/2,+\infty)) = 0$ because $s$ is a finite Borel measure on $[0,+\infty)$ with no mass point at zero (since $|\widehat{s}(0)| < \infty$). On the other hand,
\[ |\int_{\Gamma_{\varepsilon,1,\theta}}   \int_{[0,2\varepsilon) }  \frac{d s(x) }{\zeta - x}\frac{g(\zeta)d \zeta}{\zeta - z}| \leq
\int_{\Gamma_{\varepsilon,1,\theta}}   \int_{[0,2\varepsilon) } \frac{|\zeta x|}{|\zeta - x|} \frac{d| s|(x)}{|x|}  \frac{|\zeta g(\zeta)|}{|\zeta - z|}\frac{|d \zeta|}{|\zeta|^2} \leq
\]
\[ \frac{2\varepsilon}{\sin \theta} \sup_{\zeta \in \Gamma_{\varepsilon,1,\theta}, z \in {\mathcal{K}}} \frac{|\zeta g(\zeta)|}{|\zeta - z|} \int_{\varepsilon}^1 \frac{dt}{t^2} \int_{[0,2\varepsilon)  } \frac{d|s|(x)}{x} \longrightarrow 0, \qquad \varepsilon \to 0,
\]
because from f) it follows that $\sup_{\zeta \in \Gamma_{\varepsilon,1,\theta}, z \in {\mathcal{K}}} \frac{|\zeta g(\zeta)|}{|\zeta - z|}$ is uniformly bounded for $0 < \varepsilon <1 $, and $\int_{[0,2\varepsilon) } \frac{d|s|(x)}{x}$ tends to zero when $\varepsilon \to 0$  since $|\widehat{s}(0)| < \infty$ (use Lebesgue's dominated convergence theorem). We also have
\[ |\int_{\Gamma_{\varepsilon,1,\theta}}   \int_{(R/2,+\infty)}  \frac{d s(x) }{\zeta - x}\frac{g(\zeta)d \zeta}{\zeta - z}| \leq |\int_{\Gamma_{\varepsilon,1,\theta}}   \int_{(R/2,+\infty)}  \frac{d s(x) }{\zeta - x}\frac{\zeta g(\zeta)}{\zeta - z} \frac{d \zeta}{\zeta}| \leq
\]
\[ \frac{2}{R\sin \theta} \sup_{\zeta \in \Gamma_{\varepsilon,1,\theta}, z \in {\mathcal{K}}} \frac{|\zeta g(\zeta)|}{|\zeta - z|} \int_{\varepsilon}^1 \frac{dt}{t} \int_{(R/2,+\infty)} d|s|(x) \longrightarrow 0, \quad \qquad R = 1/\varepsilon \to \infty.
\]

Finally, from Fubini's theorem and Cauchy's integral formula, we have that
\[\int_{\Gamma_{\varepsilon,R,\theta}}   \int_{[2\varepsilon,R/2]}  \frac{d s(x) }{\zeta - x}\frac{g(\zeta)d \zeta}{\zeta - z}  =
\]
\[\int_{C_{\varepsilon,\theta,2} \cup \Gamma_{\varepsilon,R,\theta} \cup C_{R,\theta,2}}   \int_{[2\varepsilon,R/2]}  \frac{d s(x) }{\zeta - x}\frac{g(\zeta)d \zeta}{\zeta - z}  -
\int_{C_{\varepsilon,\theta,2} \cup C_{R,\theta,2}}   \int_{[2\varepsilon,R/2]}  \frac{d s(x) }{\zeta - x}\frac{g(\zeta)d \zeta}{\zeta - z} =
\]
\[\int_{[2\varepsilon,R/2]} \int_{C_{\varepsilon,\theta,2} \cup \Gamma_{\varepsilon,R,\theta} \cup C_{R,\theta,2}}     \frac{g(\zeta)d \zeta}{\zeta - z}  \frac{d s(x) }{\zeta - x} -
\int_{C_{\varepsilon,\theta,2} \cup C_{R,\theta,2}}   \int_{[2\varepsilon,R/2]}  \frac{d s(x) }{\zeta - x}\frac{g(\zeta)d \zeta}{\zeta - z} =
\]
\[ 2 \pi i \int_{[2\varepsilon,R/2]}  \frac{g(x) d s(x) }{z - x} -
\int_{C_{\varepsilon,\theta,2} \cup C_{R,\theta,2}}   \int_{[2\varepsilon,R/2]}  \frac{d s(x) }{\zeta - x}\frac{g(\zeta)d \zeta}{\zeta - z}.
\]
The assumptions on $g$ and $s$ together with the Lebesgue dominated convergence theorem imply that the first integral on the last line tends to $\int \frac{g(x) d s(x)}{z-x}$ when $\varepsilon \to 0, R \to \infty$. It rests to show that the other term tends to zero as $\varepsilon \to 0, R \to \infty$.

In fact, using e), we obtain
\[
| \int_{C_{R,\theta,2}}   \int_{[2\varepsilon,R/2]}  \frac{d s(x) }{\zeta - x}\frac{g(\zeta)d \zeta}{\zeta - z} | \leq {\mathcal{O}}(1)\frac{4\pi |s|([0,+\infty))}{|R-|z||} \longrightarrow 0, \quad R \to 0,
\]
and, on account of f),
\[
| \int_{C_{\varepsilon,\theta,2}}   \int_{[2\varepsilon,R/2]}  \frac{d s(x) }{\zeta - x}\frac{g(\zeta)d \zeta}{\zeta - z} | \leq \frac {2\pi o(1)}{||z|-\varepsilon|} \int_{[2\varepsilon,R/2]} \frac{d |s| (x)}{x - \varepsilon} \longrightarrow 0, \quad \varepsilon \to 0,
\]
since $x- \varepsilon \geq x/2$ for $x \geq 2\varepsilon$ and $|\widehat{s}(0)| < \infty$. We are done. \hfill $\Box$

\vspace{0,2cm}

We will be applying the previous lemma on products and quotients of Cauchy transforms so it is convenient to point out some properties of these functions.

\bl \label{lem:otro} Let $s \in {\mathcal{M}}([a,+\infty))$. Then, for each $\theta > 0$
\[ \sup_{z \in C_{R,\theta,1}} |\widehat{s}(z)| = {\mathcal{O}}(1/R), \qquad R \to \infty,
\]
and, if $a$ is  not a mass point of $s$,
\[\sup_{z \in C_{\varepsilon,\theta,1}} |(z-a)\widehat{s}(z)| = {o}(1), \qquad \varepsilon \to 0.
\]
\el

{\bf Proof.} Without loss of generality we can assume that $a=0$. Then
\[ \sup_{z \in C_{R,\theta,1}} |\int \frac{d s (x)}{z-x} | \leq \frac{|s|[0,+\infty)}{R} \sup_{z \in C_{R,\theta,1}} |\frac{z}{z- x }| \leq \frac{|s|[0,+\infty)}{R \sin \theta},
\]
giving the first relation. On the other hand,
\[   \sup_{z \in C_{\varepsilon,\theta,1}} |\int \frac{z d s (x)}{z-x} | \leq     \int \sup_{z \in C_{\varepsilon,\theta,1}} | \frac{z }{z-x} | d |s| (x) \leq     \int  | \frac{z_{\varepsilon}(x) }{z_{\varepsilon}(x)-x} | d |s| (x),
\]
where $z_{\varepsilon}(x) \in C_{\varepsilon,\theta,1}$. Since $|\frac{z_{\varepsilon}(x) }{z_{\varepsilon}(x)-x} | \leq \frac{1}{\sin \theta}, x \in C_{\varepsilon,\theta,1}$, Lebesgue's dominated convergence theorem implies that
\[ \lim_{\varepsilon \to 0} \sup_{z \in C_{\varepsilon,\theta,1}} |\int \frac{z d s (x)}{z-x} |= \int g(x) d|s|(x) =0
\]
where $g(0) =1$ and $g(x) =0, x \in (0,+\infty).$ \hfill $\Box$

\vspace{0,2cm}

\bl \label{alfabeta} Let $\sigma_{\alpha} \in
{\mathcal{M}}(\Delta_{\alpha}), \sigma_{\beta} \in
{\mathcal{M}}(\Delta_{\beta})$, and $\la \sigma_{\alpha}, {\sigma}_{\beta} \ra\in {\mathcal{M}}(\Delta_{\alpha})$.
Then,  $\la \sigma_{\beta}, {\sigma}_{\alpha} \ra\in {\mathcal{M}}(\Delta_{\beta})$ and
\begin{equation} \label{2.1} \widehat{\sigma}_{\alpha} (z)\widehat{\sigma}_{\beta}(z)
=\langle \sigma_{\alpha},\sigma_{\beta} \widehat{\rangle}(z)+\langle
\sigma_{\beta},{\sigma}_{\alpha}\widehat{\ra}(z), \quad z \in
{\mathbb{C}} \setminus \left(\Delta_{\alpha} \cup
\Delta_{\beta}\right),
\end{equation}
\el

{\bf Proof.}
In fact,
\[
\widehat{\sigma}_{\alpha} (z)\widehat{\sigma}_{\beta} (z)=\int \int
\frac{d \sigma_{\alpha}(x_{\alpha}) d
\sigma_{\beta}(x_{\beta})}{(z-x_{\alpha})(z-x_{\beta})}=\int \int
\left(\frac{1}{z-x_{\alpha}}-\frac{1}{z-x_{\beta}}\right)\frac{d
\sigma_{\alpha}(x_{\alpha}) d
\sigma_{\beta}(x_{\beta})}{x_{\alpha}-x_{\beta}}.
\]
Therefore, the right hand is finite because $\widehat{\sigma}_{\alpha} \in
{\mathcal{M}}(\Delta_{\alpha})$ and $\widehat{\sigma}_{\beta} \in
{\mathcal{M}}(\Delta_{\beta})$. We also have that
\[ \int \int \frac{d \sigma_{\alpha}(x_{\alpha}) d
\sigma_{\beta}(x_{\beta})}{(z-x_{\alpha})(x_{\alpha}-x_{\beta})} = \int \frac{\widehat{\sigma}_{\beta}(x_{\alpha})d \sigma_{\alpha}(x_{\alpha})}{z-x_{\alpha}}
\]
is finite, since $\widehat{\sigma}_{\beta} \in L_1(\sigma_{\alpha})$. Consequently,
\[  \int \int
\left(\frac{1}{z-x_{\alpha}}-\frac{1}{z-x_{\beta}}\right)\frac{d
\sigma_{\alpha}(x_{\alpha}) d
\sigma_{\beta}(x_{\beta})}{x_{\alpha}-x_{\beta}}
\]
can be separated in two integrals obtaining (\ref{2.1}). From (\ref{2.1}) it follows that
\[ |\la \sigma_{\alpha}, \sigma_{\beta} {\ra}| = \lim_{z \to \infty} z \la \sigma_{\alpha}, \sigma_{\beta} \widehat{\ra}(z) = \lim_{z \to \infty} z \widehat{\sigma}_{\alpha}(z) \widehat{\sigma}_{\beta}(z) - \lim_{z \to \infty} z \la \sigma_{\beta}, \sigma_{\alpha} \widehat{\ra}(z) = - |\la \sigma_{\beta}, \sigma_{\alpha} {\ra}|
\]
is finite and we conclude the proof ($|s|$ denotes the total mass of the measure $s$).
\hfill $\Box$

\vspace{0,2cm}

\bl Let $\sigma_{\alpha} \in
{\mathcal{M}}(\Delta_{\alpha}), \sigma_{\beta} \in
{\mathcal{M}}(\Delta_{\beta})$, $\Delta_{\alpha} \cap \Delta_{\beta} = \{x_{\alpha,\beta}\}$ and $\widehat{\sigma}_{\beta}(x_{\alpha,\beta}) \in {\mathbb{R}}$. Then, $\la \sigma_{\alpha}, {\sigma}_{\beta} \ra\in {\mathcal{M}}(\Delta_{\alpha})$  and $\la \sigma_{\beta}, {\sigma}_{\alpha} \ra\in {\mathcal{M}}(\Delta_{\beta})$ .
\el

{\bf Proof.} By monotonicity, we have that
\[ |\la \sigma_{\alpha},\sigma_{\beta} \ra| = |\int \widehat{\sigma}_{\beta}(x_{\alpha}) d \sigma_{\alpha}(x_{\alpha})| =
\int |\widehat{\sigma}_{\beta}(x_{\alpha})| d |\sigma_{\alpha}|(x_{\alpha}) \leq |\widehat{\sigma}_{\beta}(x_{\alpha,\beta})| |\sigma_{\alpha}|
\]
is finite. Thus, $\la \sigma_{\alpha}, {\sigma}_{\beta} \ra\in {\mathcal{M}}(\Delta_{\alpha})$  and the previous lemma implies that $\la \sigma_{\beta}, {\sigma}_{\alpha} \ra\in {\mathcal{M}}(\Delta_{\beta})$.  \hfill $\Box$

\vspace{0,2cm}

The following result is similar to \cite[Lemma 3.1]{LF4}.

\bl \label{alfabeta+} Suppose that $\Delta_{\alpha},\Delta_{\beta},$ are two intervals which have at most a common end point. Let $\sigma_{\alpha} \in
{\mathcal{M}}(\Delta_{\alpha}), \sigma_{\beta} \in
{\mathcal{M}}(\Delta_{\beta})$, and $\la \sigma_{\alpha}, {\sigma}_{\beta} \ra\in {\mathcal{M}}(\Delta_{\alpha})$.
Assume that the first three moments of $\sigma_{\alpha}$ are finite. If $\Delta_{\alpha} \cap \Delta_{\beta} = \{x_{\alpha,\beta}\}$, the point $x_{\alpha,\beta}$ is not a mass point of $\sigma_{\alpha}$ or $\sigma_{\beta}$.
Then, $\la \sigma_{\beta},\sigma_{\alpha}\widehat{\ra}/\widehat{\sigma}_{\beta} \in L_1(\tau_{\alpha,\beta})$, $\la \sigma_{\beta},\sigma_{\alpha}\widehat{\ra}  \in L_1(\tau_{\alpha,\alpha})$, and
\begin{equation} \label{2.2}
\frac{\widehat{\sigma}_{\alpha}(z)}{\langle
\sigma_{\alpha},\sigma_{\beta}\widehat{\rangle}(z)}=
\frac{|\sigma_{\alpha}|}{|\langle\sigma_{\alpha}\sigma_{\beta}\rangle|}+
\int \frac{\langle
\sigma_{\beta},{\sigma}_{\alpha}\widehat{\ra}(x_{\alpha})}{
\widehat{\sigma}_{\beta} (x_{\alpha})}\frac{d \tau_{\alpha,\beta}
(x_{\alpha})}{z-x_{\alpha}} =
\frac{|\sigma_{\alpha}|}{|\langle\sigma_{\alpha},\sigma_{\beta}\rangle|}+\langle
 \frac{{\tau}_{\alpha,\beta}}{\widehat{\sigma}_{\beta}},\sigma_{\beta}, {\sigma}_{\alpha}\widehat{\ra}(z),
\end{equation}
\begin{equation} \label{2.3}
\frac{\langle
\sigma_{\alpha},\sigma_{\beta}\widehat{\rangle}(z)}{\widehat
{\sigma}_{\alpha} (z)}=\frac{|\langle
\sigma_{\alpha},\sigma_{\beta}\rangle|}{|\sigma_{\alpha}|}- \int
\frac{{\langle \sigma_{\beta},
{\sigma}_{\alpha}\widehat{\ra}(x_{\alpha})}d \tau_{\alpha,\alpha}
(x_{\alpha})}{z-x_{\alpha}}=\frac{|\langle
\sigma_{\alpha},\sigma_{\beta}\rangle|}{|\sigma_{\alpha}|} -\langle
 {\tau}_{\alpha,\alpha},\sigma_{\beta}, {\sigma}_{\alpha}\widehat{\ra}(z).
\end{equation}
Consequently,
\begin{equation} \label{doslim}
\lim_{x \to x_{\alpha,\beta}, x \in {\mathbb{R}}\setminus \Delta_{\alpha}}\frac{\widehat{\sigma}_{\alpha}(x)}{\langle
\sigma_{\alpha},\sigma_{\beta}\widehat{\rangle}(x)} \in {\mathbb{R}}\qquad \mbox{and} \qquad \lim_{x \to \infty, x \in {\mathbb{R}} \setminus \Delta_{\alpha}} \frac{\langle
\sigma_{\alpha},\sigma_{\beta}\widehat{\rangle}(x)}{\widehat
{\sigma}_{\alpha} (x)}\in {\mathbb{R}}
\end{equation}.
\el

{\bf Proof.} For the proof of (\ref{2.2})-(\ref{2.3}) we restrict our attention to the most complicated case when $\Delta_{\alpha}$ and $\Delta_{\beta}$ are both unbounded and have a common end point. Without loss of generality we can assume that $\Delta_{\alpha}= {\mathbb{R}}_+ = [0,+\infty), \Delta_{\beta}= {\mathbb{R}}_- = (-\infty,0]$ and thus $x_{\alpha.\beta} =0$. Set
\[
w(z) = \frac{\widehat{\sigma}_{\alpha}(z)}{\langle
\sigma_{\alpha},\sigma_{\beta}\widehat{\rangle}(z)}-
\frac{|\sigma_{\alpha}|}{|\langle\sigma_{\alpha}\sigma_{\beta}\rangle|}
\in {\mathcal{H}} \left({\mathbb{C}} \setminus
{\mathbb{R}}_+\right).
\]
From (\ref{eq:1}) and (\ref{2.1}) (see also the last statement of Lemma \ref{momentos} and Remark \ref{rem:1}),   there exists $\tau_{\alpha,\beta} \in {\mathcal{M}}(\Delta_{\alpha})$ such that
\[
\frac{\widehat{\sigma}_{\alpha}(z)}{\langle
\sigma_{\alpha},\sigma_{\beta}\widehat{\rangle}(z)} = \frac{\widehat
{\sigma}_{\beta} (z) \widehat{\sigma}_{\alpha} (z)}{\widehat{
\sigma}_{\beta}(z)\langle
\sigma_{\alpha},\sigma_{\beta}\widehat{\rangle}(z)}  =\frac{\langle
\sigma_{\alpha},\sigma_{\beta}\widehat{\rangle}(z)+\langle
\sigma_{\beta}, {\sigma}_{\alpha}\widehat{\ra}(z)}{\widehat
{\sigma}_{\beta} (z)\langle
\sigma_{\alpha},\sigma_{\beta}\widehat{\rangle}(z)} =
\]
\[
\frac{1}{\widehat {\sigma}_{\beta} (z)} +  \frac{\langle
\sigma_{\beta}, {\sigma}_{\alpha}\widehat{\ra}(z)}{\widehat
{\sigma}_{\beta} (z)}\ell_{\alpha,\beta} + \frac{\langle
\sigma_{\beta}, {\sigma}_{\alpha}\widehat{\ra}(z)}{\widehat
{\sigma}_{\beta} (z)} \widehat{\tau}_{\alpha,\beta} (z).
\]
Consequently,
\[ w(z) = h_1(z) + g(z)\widehat{\tau}_{\alpha,\beta} (z), \qquad z \in {\mathbb{C}} \setminus {\mathbb{R}},
\]
where
\[ h_1(z) = -\frac{|\sigma_{\alpha}|}{|\langle\sigma_{\alpha}\sigma_{\beta}\rangle|} + \frac{1}{\widehat {\sigma}_{\beta} (z)} +  \frac{\langle
\sigma_{\beta}, {\sigma}_{\alpha}\widehat{\ra}(z)}{\widehat
{\sigma}_{\beta} (z)}\ell_{\alpha,\beta}, \qquad g(z) = \frac{\langle
\sigma_{\beta}, {\sigma}_{\alpha}\widehat{\ra}(z)}{\widehat
{\sigma}_{\beta} (z)}
\]
$h_1 \in {\mathcal{H}}({\mathbb{C}} \setminus {\mathbb{R}}_-), g \in {\mathcal{H}}({\mathbb{C}} \setminus {\mathbb{R}}_-)$.

Given a compact set ${\mathcal{K}} \subset {\mathbb{C}} \setminus [0,+\infty)$, fix $\theta >0$ sufficiently small so that ${\mathcal{K}}$ lies to the left of $\Gamma_{\theta}$. The independent term of the asymptotic expansion at $\infty$ of $ {\widehat{\sigma}_{\alpha}(z)}/{\langle
\sigma_{\alpha},\sigma_{\beta}\widehat{\rangle} }$ equals ${|\sigma_{\alpha}|}/{|\langle\sigma_{\alpha}\sigma_{\beta}\rangle|}$. It readily follows that  $w$ satisfies c) of Lemma \ref{lem:2}. On the other hand, using the second part of Lemma \ref{lem:otro} we have that $w$ also verifies  d) of Lemma \ref{lem:2}. By (\ref{eq:2.2})
\[ w(z) = \frac{1}{2\pi i} \int_{\Gamma_{\theta}} \frac{w(\zeta) d \zeta}{\zeta - z} = \frac{1}{2\pi i} \int_{\Gamma_{\theta}} \frac{(h_1(\zeta) + g(\zeta) \widehat{\tau}_{\alpha,\beta}(\zeta))d \zeta}{\zeta - z} .
\]
For each $z \in {\mathcal}{K},$ the function $h(\zeta) = h_1(\zeta)/(\zeta -z)$ as a function of $\zeta$
satisfies a) and b) of Lemma \ref{lem:2}, $ g(z)$ satisfies e) and f), and $|\widehat{\tau}_{\alpha,\beta}(0)| < \infty$. So, (\ref{2.2}) holds true on account of (\ref{eq:2.1}) and (\ref{eq:2.3}).

Because of the right hand in (\ref{2.2}), we have that $\widehat{\sigma}_{\alpha}(t)/\la \sigma_{\alpha},\sigma_{\beta} \widehat{\ra}(t)$ is monotonic for $t \in {\mathbb{R}}$ as $t$ approaches $0$. Therefore, $\lim_{t \to 0-} \widehat{\sigma}_{\alpha}(t)/\la \sigma_{\alpha},\sigma_{\beta} \widehat{\ra}(t)$ exists if the function is bounded above in absolute value on the negative real axis when $t$ is sufficiently close to zero. Fix $A> 0$. There exist positive constants $C_1, C_2,C_3,$ such that $|\widehat{\sigma}_{\beta}(x)| > C_1, x \in [0,A],$ $\int_A^{+\infty} \frac{d|\sigma_{\alpha}|(x)}{|t -x|} \leq \int_A^{+\infty} \frac{|d\sigma_{\alpha}|(x)}{x} \leq C_2,$ and $\int_0^A \frac{d|\sigma_{\alpha}|(x)}{|t -x|} \geq C_3, t \in [-1,0)$; consequently,
\[ \frac{|\widehat{\sigma}_{\alpha}(t)|}{|\la \sigma_{\alpha},\sigma_{\beta} \widehat{\ra}(t)|} = \frac{\int \frac{d|\sigma_{\alpha}|(x)}{|t - x|}}{\int \frac{|\widehat{\sigma}_{\beta}(x)|d|\sigma_{\alpha}|(x)}{|t - x|}} \leq \frac{\int_0^A \frac{d|\sigma_{\alpha}|(x)}{|t - x|} + C_2}{C_1 \int_0^A \frac{d|\sigma_{\alpha}|(x)}{|t - x|} } \leq \frac{1}{C_1} + \frac{C_2}{C_1C_3}, \qquad t \in [-1,0),
\]
as we needed to prove. This settles the first part of the last statement.

The proof of (\ref{2.3}) is similar to that of (\ref{2.2}). It is based on the fact that
\[
\frac{\langle\sigma_{\alpha},\sigma_{\beta}\widehat{\rangle}(z)}{
\widehat{\sigma}_\alpha(z)}-\frac{|\langle
\sigma_{\alpha},\sigma_{\beta}\rangle|}{|\sigma_{\alpha}|}
\in {\mathcal{H }}\left({\mathbb{C}} \setminus
{\mathbb{R}}_+ \right),
\]
and that (\ref{eq:1}) and (\ref{2.1}) imply the formula
\[
\frac{\langle
\sigma_{\alpha},\sigma_{\beta}\widehat{\rangle}(z)}{\widehat{
\sigma}_{\alpha}(z)} =\frac{\widehat{\sigma}_{\alpha}
(z)\widehat{\sigma}_{\beta} (z)-\langle \sigma_{\beta},
{\sigma}_{\alpha}\widehat{\ra}(z)}{\widehat{ \sigma}_{\alpha}(z)} =
\widehat{\sigma}_{\beta} (z) - \langle \sigma_{\beta},
{\sigma}_{\alpha}\widehat{\ra}(z)\ell_{\alpha,\alpha}(z) - \langle
\sigma_{\beta},
{\sigma}_{\alpha}\widehat{\ra}(z)\widehat{\tau}_{\alpha,\alpha}(z).
\]
Take $w(z) = \frac{\langle
\sigma_{\alpha},\sigma_{\beta}\widehat{\rangle}(z)}{\widehat{
\sigma}_{\alpha}(z)} -\frac{|\langle
\sigma_{\alpha},\sigma_{\beta}\rangle|}{|\sigma_{\alpha}|}, h(z) = \widehat{\sigma}_{\beta} (z) - \langle \sigma_{\beta},
{\sigma}_{\alpha}\widehat{\ra}(z)\ell_{\alpha,\alpha}(z) -\frac{|\langle
\sigma_{\alpha},\sigma_{\beta}\rangle|}{|\sigma_{\alpha}|}, g(z) = - \langle
\sigma_{\beta},
{\sigma}_{\alpha}\widehat{\ra}(z),$ $s = \tau_{\alpha,\alpha}$, and proceed as in the previous case.

The second limit in (\ref{doslim}) is derived using arguments similar to those employed in proving the first limit. We leave it to the reader.
\hfill $\Box$ \vspace{0,2cm}

Let us extend Lemma \ref{alfabeta}. The result will not be employed but gives an interesting relation.

\bl Let $(s_{1,1},\ldots,s_{1,m}) = {\mathcal{N}}(\sigma_1,\ldots,\sigma_m), m \geq 2,$ be given. Then, all the measures appearing in the formula below are finite  and for all $z \in {\mathbb{C}} \setminus (\Delta_1 \cup \Delta_m)$
\begin{equation} \label{extendido}
\la \sigma_m,\ldots,\sigma_1 \widehat{\ra}(z) - \sum_{k=1}^{m-1} (-1)^{k}\la \sigma_m,\ldots,\sigma_{k+1}\widehat{\ra}(z)\la \sigma_1,\ldots,\sigma_{k}\widehat{\ra}(z) + (-1)^{m}\la \sigma_1,\ldots,\sigma_{m}\widehat{\ra}(z)\equiv 0.
\end{equation}
In particular,  $|\la \sigma_m,\ldots,\sigma_1  {\ra}| = (-1)^{m-1} |\la \sigma_1,\ldots,\sigma_{m} {\ra}| $. Moreover, $(\sigma_m,\ldots,\sigma_1)$ is a generator of a Nikishin system.
\el

{\bf Proof.} This result reduces to Lemma \ref{alfabeta} when $m=2$. The general case follows by induction.
Let us assume that the lemma is true for any Nikishin system with $m-1, m \geq 3,$ generating measures and let
us show that it also holds when the number of generating measures is $m$.

We have
\[ \frac{1}{(x_2 -x_1)\cdots (x_m - x_{m-1})(z - x_m)} - \frac{(-1)^{m-1}}{(z - x_1)(x_1 -x_2)\cdots (x_{m-2} - x_{m-1})(x_{m-1} - x_m)} =
\]
\[ \frac{x_m - x_1}{(z - x_1)(x_2 -x_1)\cdots (x_m - x_{m-1})(z - x_m)}  = \sum_{k =1}^{m-1}  \frac{x_{k+1} - x_{k}}{(z - x_1)(x_2 -x_1)\cdots (x_m - x_{m-1})(z - x_m)}.
\]
For $z \in {\mathbb{C}} \setminus \Delta_1 $
\[ \int \cdots \int \frac{d\sigma_1(x_1)\cdots d\sigma_m(x_m) }{(z - x_1)(x_1 -x_2)\cdots (x_{m-2} - x_{m-1})(x_{m-1} - x_m)} = \la \sigma_1,\ldots,\sigma_m \widehat{\ra} (z)
\]
is finite. By hypothesis and induction hypothesis, for  $k=1,\ldots,m-1,$ and $z \in {\mathbb{C}} \setminus (\Delta_1 \cup \Delta_m)$
\[ \int \cdots \int \frac{(x_{k+1} - x_{k})d\sigma_1(x_1)\cdots d\sigma_m(x_m) }{(z - x_1)(x_2 -x_1)\cdots (x_m - x_{m-1})(z - x_m)} = (-1)^{k-1} \la \sigma_m ,\ldots, \sigma_{k+1} \widehat{\ra} (z) \la \sigma_1,\ldots,\sigma_k \widehat{\ra}
\]
is also finite. Integrating term by term, it follows that for $z \in {\mathbb{C}} \setminus (\Delta_1 \cup \Delta_m)$
\[ \int \cdots \int \frac{d\sigma_1(x_1)\cdots d\sigma_m(x_m) }{(x_2 -x_1)\cdots (x_m - x_{m-1})(z - x_m)} = \la \sigma_m,\ldots,\sigma_1 \widehat{\ra}(z)
\]
is finite and satisfies (\ref{extendido}).

Using (\ref{extendido}) and Lebesgue's dominated convergence theorem, we obtain that
\[ |\la \sigma_m,\ldots,\sigma_1  {\ra}| = \lim_{t \to \infty} it \la \sigma_m,\ldots,\sigma_1  \widehat{\ra}(it) =
\]
\[(-1)^{m-1} \lim_{t \to \infty} it \la \sigma_1,\ldots,\sigma_{m} \widehat{\ra}(it) = (-1)^{m-1} |\la \sigma_1,\ldots,\sigma_{m} {\ra}|.
\]
Consequently, $\la \sigma_m,\ldots,\sigma_1  {\ra} \in {\mathcal{M}}(\Delta_m)$ and the relation between the total variation of $\la \sigma_m,\ldots,\sigma_1  {\ra}$ and $\la \sigma_1,\ldots,\sigma_m  {\ra}$ has been established.

We have proved that  $\la \sigma_k, \sigma_{k-1}, \ldots,\sigma_{j}\ra \in {\mathcal{M}}(\Delta_k) $ for all $1 \leq j \leq k \leq m$. The assumption that common end points of consecutive intervals $\Delta_j$ are not mass points  of the corresponding measures remains valid. Therefore, $(\sigma_m,\ldots,\sigma_1)$ is a generator of a Nikishin system. \hfill $\Box$

\vspace{0,2cm}
An iterated application of (\ref{2.2})-(\ref{2.3}) allows to derive important formulas which we present in the next result analogous to  \cite[Lemma 3.2]{LF4}.

\bl \label{cocientes} Let $(s_{1,1},\ldots,s_{1,m}) =
{\mathcal{N}}(\sigma_1,\ldots,\sigma_m)$ be given. Assume that  the first three moments of all the generating measures are finite. Then, all the measures appearing in the formulas below are finite, and:
\begin{equation} \label{4.4}
\frac{\widehat{s}_{1,k}}{\widehat{s}_{1,1}} =
\frac{|s_{1,k}|}{|s_{1,1}|} - \la \tau_{1,1},\la s_{2,k},\sigma_1
\ra \widehat{\ra}  , \qquad  1=j < k \leq m,
\end{equation}
\begin{equation} \label{4.5} \frac{\widehat{s}_{1,k}}{\widehat{s}_{1,j}} =
\frac{|s_{1,k}|}{|s_{1,j}|} + (-1)^j \la
\tau_{1,j},\la\tau_{2,j},s_{1,j} \ra,\ldots,
\la\tau_{j,j},s_{j-1,j}\ra,\la s_{j+1,k},\sigma_j\ra \widehat{\ra},
\quad 2\leq j < k \leq m,
\end{equation}
\begin{equation} \label{4.6}
\frac{\widehat{s}_{1,1}}{\widehat{s}_{1,j}} =
\frac{|s_{1,1}|}{|s_{1,j}|} + \la
\frac{\tau_{1,j}}{\widehat{s}_{2,j}}, \la s_{2,j},\sigma_1\ra
\widehat{\ra}=
\end{equation}
\[ \frac{|s_{1,1}|}{|s_{1,j}|} + \frac{|\la
s_{2,j},\sigma_1\ra|}{|s_{2,j}|} \widehat{\tau}_{1,j} - \la
\tau_{1,j},\la \tau_{2,j},s_{1,j}\ra \widehat{\ra}, \qquad 1 =k <
j\leq m,
\]
\begin{equation} \label{4.7}
\frac{\widehat{s}_{1,2}}{\widehat{s}_{1,j}} =
\frac{|s_{1,2}|}{|s_{1,j}|} - \la \tau_{1,j}, \frac{\la
\tau_{2,j},s_{1,j}\ra}{\widehat{s}_{3,j}}, \la s_{3,j},\sigma_2\ra
\widehat{\ra}=
\end{equation}
\[ \frac{|s_{1,2}|}{|s_{1,j}|} -
\frac{|\la s_{3,j},\sigma_2\ra|}{|s_{3,j}|}  \la {\tau}_{1,j},\la
\tau_{2,j},s_{1,j}\ra \widehat{\ra} + \la
\tau_{1,j},\la \tau_{2,j},s_{1,j}\ra , \la \tau_{3,j},s_{2,j}\ra
\widehat{\ra},\qquad 2 =k < j \leq m,
\]
\begin{equation} \label{4.2}
\frac{\widehat{s}_{1,k}}{\widehat{s}_{1,j}} =
\frac{|s_{1,k}|}{|s_{1,j}|} + (-1)^{k-1} \la
\tau_{1,j},\la\tau_{2,j},s_{1,j} \ra,\ldots,
\la\tau_{k-1,j},s_{k-2,j}\ra, \frac{\la \tau_{k,j},
{s}_{k-1,j}\ra}{\widehat{s}_{k+1,j}},\la s_{k+1,j},\sigma_k\ra
\widehat{\ra} =
\end{equation}
\[ \frac{|s_{1,k}|}{|s_{1,j}|} + (-1)^{k-1} \frac{|\la s_{k+1,j},\sigma_k\ra|}{|s_{k+1,j}|}\la
\tau_{1,j},\la\tau_{2,j},s_{1,j} \ra,\ldots,
\la\tau_{k-1,j},s_{k-2,j}\ra, \la\tau_{k,j},s_{k-1,j}\ra
\widehat{\ra} \,+
\]
\[ (-1)^{k}  \la \tau_{1,j},\la\tau_{2,j},s_{1,j} \ra,\ldots,
 \la\tau_{k,j},s_{k-1,j}\ra, \la\tau_{k+1,j},s_{k,j}\ra
\widehat{\ra}\,. \qquad 3=k < j \leq m.
\]
\el

{\bf Proof.}  The main difference in the proof of this lemma with respect to that of \cite[Lemma 3.2]{LF4} is that one must verify that the new assumptions allow to apply formulas (\ref{2.2})-(\ref{2.3})  as many times as required. We begin pointing out that Lemmas \ref{alfabeta} and \ref{alfabeta+} guarantee that the measures between commas in formulas (\ref{4.4})-(\ref{4.2}) are finite. That their products, as indicated, give rise to finite measures is proved step by step.

The independent term of the asymptotic expansion at $\infty$ of the ratios of Cauchy transforms appearing on the left of (\ref{4.4})-(\ref{4.2}) is the constant appearing on the right. If we prove that those ratios of Cauchy transforms may be expressed as a constant plus a Cauchy transform of a finite measure (whose asymptotic expansion at infinity have independent term equal to zero) we have that the constant has to be the one given. Consequently, we will not
pay attention to the constants coming out of the consecutive
transformations that we make in our deduction and simply denote them with
consecutive constants $C_j$.

Obviously,  (\ref{4.4}) is deduced from (\ref{2.3}) taking
$\sigma_{\alpha} = \sigma_1 = s_{1,1} $ and $\sigma_{\beta} = \la
\sigma_2,\cdots,\sigma_k \ra = s_{2,k}$. Formula (\ref{4.5}) is
obtained applying (\ref{2.3}) inside out several times  as we show.

Let $2 \leq j < k \leq m$. Using (\ref{2.3}) on $
{\widehat{s}_{j,k}}/{\widehat{s}_{j,j}}$, we have that
\begin{equation} \label{4.3}
\la \sigma_{j-1},\sigma_j,\ldots,\sigma_k\widehat{\ra} = \la
\frac{s_{j-1,j}}{\widehat{s}_{j,j}}, s_{j,k}\widehat{\ra} =   \la
\frac{\widehat{s}_{j,k}}{\widehat{s}_{j,j}}s_{j-1,j}\widehat{\ra} =
C_1 \widehat{s}_{j-1,j} - \la s_{j-1,j},\tau_{j,j}, \la
s_{j+1,k},\sigma_j\ra \widehat{\ra}.
\end{equation}
The last equality is due to the fact that $\la \sigma_{j-1},\sigma_j,\ldots,\sigma_k\widehat{\ra}$ and $\widehat{s}_{j-1,j}$ are finite so we can integrate term by term. Multiply (\ref{4.3}) by $z$ and make $z$ tend to infinity conveniently, to deduce that $\la s_{j-1,j},\tau_{j,j}, \la
s_{j+1,k},\sigma_j\ra  {\ra}$ is finite.  In particular, if $j=2$ we get
\[ \la \sigma_{1},\sigma_2,\ldots,\sigma_k\widehat{\ra} = \la \frac{s_{1,2}}{\widehat{s}_{2,2}},
s_{2,k}\widehat{\ra} =   \la
\frac{\widehat{s}_{2,k}}{\widehat{s}_{2,2}}s_{1,2}\widehat{\ra} =
C_1 \widehat{s}_{1,2} - \la s_{1,2},\tau_{2,2}, \la
s_{3,k},\sigma_2\ra \widehat{\ra},
\]
and applying (\ref{2.3}) on ${\la s_{1,2},\tau_{2,2}, \la
s_{3,k},\sigma_2\ra \widehat{\ra}}/{\widehat{s}_{1,2}}$, it follows
that
\[ \frac{\widehat{s}_{1,k}}{\widehat{s}_{1,2}} = C_1 - \frac{1}{\widehat{s}_{1,2}} \la s_{1,2},\tau_{2,2}, \la
s_{3,k},\sigma_2\ra \widehat{\ra} = \frac{|s_{1,k}|}{|s_{1,2}|} +
\la \tau_{1,2},\la \tau_{2,2}, s_{1,2} \ra, \la s_{3,k},\sigma_2\ra
\widehat{\ra}
\]
which is (\ref{4.5}) for $j=2$.

Assume that $j \geq 3$. Using (\ref{4.3}), we can write
\[ \widehat{s}_{j-2,k} = \la \sigma_{j-2},\frac{\widehat{s}_{j,k}}{\widehat{s}_{j,j}}s_{j-1,j}\widehat{\ra} =  C_1  \widehat{s}_{j-2,j} - \la \sigma_{j-2}, s_{j-1,j}, \tau_{j,j}, \la s_{j+1,j},\sigma_j \ra \widehat{\ra} =
\]
\[ C_1  \widehat{s}_{j-2,j} - \la \frac{s_{j-2,j}}{\widehat{s}_{j-1,j}}, s_{j-1,j}, \tau_{j,j}, \la s_{j+1,j},\sigma_j \ra \widehat{\ra}.
\]
Applying (\ref{2.3}) on $\la   s_{j-1,j}, \tau_{j,j}, \la s_{j+1,j},\sigma_j \ra \widehat{\ra}/{\widehat{s}_{j-1,j}}$ and integrating term by term, it follows
\begin{equation} \label{a}
\widehat{s}_{j-2,k} = C_2 \widehat{s}_{j-2,j} +  \la s_{j-2,j}, \tau_{j-1,j}, \la \tau_{j,j},s_{j-1,j}\ra,\la
s_{j+1,j},\sigma_j \ra \widehat{\ra}.
\end{equation}
Term by term integration is legal because ${s}_{j-2,k}$ and ${s}_{j-2,j}$ are finite. Multiplying (\ref{a}) by $z$ and making $z$ tend to infinity conveniently, we obtain that $\la s_{j-2,j}, \tau_{j-1,j}, \la \tau_{j,j},s_{j-1,j}\ra,\la
s_{j+1,j},\sigma_j \ra \widehat{\ra}$ is finite. If $j=3$, (\ref{a}) reduces to
\[\frac{\widehat{s}_{1,k}}{\widehat{s}_{1,3}} = C_2  +  \la s_{1,3}, \tau_{2,3}, \la \tau_{3,3},s_{2,3}\ra,\la
s_{4,3},\sigma_3 \ra \widehat{\ra}/\widehat{s}_{1,3}
\]
and making use of (\ref{2.3}) on the second term on the right hand side we get (\ref{4.5}) for $j=3$. For an arbitrary $j, 1\leq j < k \leq m$, after $j-1$ steps we arrive to (\ref{4.5}).

Regarding (\ref{4.6})-(\ref{4.2}), assuming that the first equalities are true, the  second
equalities follow directly from (\ref{2.3}) since we get
\begin{equation} \label{4.8}
\frac{\la  {s}_{k+1,j},\sigma_k \widehat{\ra}}{\widehat{s}_{k+1,j}}
= \frac{|\la
 {s}_{k+1,j},\sigma_k  \ra |}{|{s}_{k+1,j}|} - \la \tau_{k+1,j},
 s_{k,j} \widehat{\ra}.
\end{equation}
When $k=1$ formula (\ref{4.6}) follows
directly from (\ref{2.2}) taking $\sigma_{\alpha}= \sigma_1 =
s_{1,1}$ and $\sigma_{\beta} = s_{2,j} $.

In order to prove the first equalities in (\ref{4.7})-(\ref{4.2}), assume that $2 \leq k < j \leq m$. Using
(\ref{2.2}) and that $\la \sigma_{k-1},\sigma_k \ra, s_{k-1,j},$ are finite, it follows that
\begin{equation} \label{b} \la \sigma_{k-1}, \sigma_k\widehat{\ra} = \la \frac{s_{k-1,j}}{\widehat{s}_{k,j}},
s_{k,k}\widehat{\ra} =   \la
\frac{\widehat{s}_{k,k}}{\widehat{s}_{k,j}}s_{k-1,j}\widehat{\ra} =
C_5 \widehat{s}_{k-1,j} + \la
s_{k-1,j},\frac{\tau_{k,j}}{\widehat{s}_{k+1,j}}, \la
s_{k+1,j},\sigma_{k} \ra \widehat{\ra}.
\end{equation}
Multiplying $(\ref{b})$ by $z$ and making $z \to \infty$ we conclude that
$\la
s_{k-1,j},\frac{\tau_{k,j}}{\widehat{s}_{k+1,j}}, \la
s_{k+1,j},\sigma_{k} \ra {\ra}$ is finite. Consequently, when $k=2$
\[ \frac{\widehat{s}_{1,2}}{\widehat{s}_{1,j}} =  C_5 + \frac{1}{\widehat{s}_{1,j}} \la
s_{1,j},\frac{\tau_{2,j}}{\widehat{s}_{3,j}}, \la s_{3,j},\sigma_2
\ra \widehat{\ra},
\]
and applying (\ref{2.3}) we obtain (\ref{4.7}). To complete the proof of (\ref{4.2}), starting out from (\ref{b})  proceed using (\ref{2.3}) repeatedly as
in proving (\ref{4.5}).  For more details, see the proof of \cite[Lemma 3.2]{LF4}. \hfill $\Box$

Let ${\tau}_{\alpha,\beta;\gamma,\gamma}$ denote the inverse measure of
$\langle \langle \sigma_{\alpha}, \sigma_{\beta} \rangle,
\sigma_{\gamma} \ra$. That is,
\[ {1}/{\langle \langle
\sigma_{\alpha}, \sigma_{\beta} \rangle,
\sigma_{\gamma}\widehat{\rangle}(z)} =
\ell_{\alpha,\beta;\gamma,\gamma}(z) +
\widehat{\tau}_{\alpha,\beta;\gamma,\gamma}(z)
\]
where $\ell_{\alpha,\beta;\gamma,\gamma}$ denotes a first degree
polynomial. This notation seems unnecessarily
complicated. It is consistent with the one used later for more general
inverse measures which will be needed. The next result is analogous to  \cite[Lemma 3.3]{LF4}.

\bl \label{freedom} Suppose that $\Delta_{\gamma},$ $\Delta_{\alpha},$ and
$\Delta_{\beta}$ are three intervals such that $\Delta_{\gamma}, \Delta_{\alpha},$ and   $\Delta_{\beta}, \Delta_{\alpha},$ have at most a common end point. Assume also that the measures intervening in the left hand sides of the two subsequent formulas are finite, where $f$ denotes a measurable function, and the first three moments of $\sigma_{\alpha}$ are finite. If $\Delta_{\alpha} \cap \Delta_{\beta} = \{x_{\alpha,\beta}\}$, the point $x_{\alpha,\beta}$ is not a mass point of $\sigma_{\alpha}$ or $\sigma_{\beta}$; likewise if $\Delta_{\alpha} \cap \Delta_{\gamma} = \{x_{\alpha,\gamma}\}$. Then, the measures on the right hand  are also finite and
\begin{equation}\label{inversa2*}
\frac{\widehat{\sigma}_{\alpha}(z)}{\langle
\sigma_{\alpha},\sigma_{\beta}\widehat{\rangle}(z)} \langle \langle
\tau_{\alpha,\alpha} ,\sigma_{\beta},\sigma_{\alpha}\rangle, f
\sigma_{\gamma}, \sigma_{\alpha} \widehat{\rangle}(z)= \la \frac{\la
\sigma_{\beta},\sigma_{\alpha}
\widehat{\ra}}{\widehat{\sigma}_{\beta}} \tau_{\alpha,\beta},
f\sigma_{\gamma}, \sigma_{\alpha},\sigma_{\beta} \widehat{\ra}(z),
\end{equation}
\begin{equation}\label{shirenu}
\frac{\langle
\sigma_{\alpha},\sigma_{\beta}\widehat{\rangle}(z)}{\langle \langle
\sigma_{\alpha}, \sigma_{\beta} \rangle,
\sigma_{\gamma}\widehat{\rangle}(z)} \la \frac{\la
\sigma_{\beta},\sigma_{\alpha}
\widehat{\ra}}{\widehat{\sigma}_{\beta}} \tau_{\alpha,\beta},
\sigma_{\gamma}, \sigma_{\alpha},\sigma_{\beta} \widehat{\ra}(z)=
\la  \frac{\langle \sigma_{\beta}
,\sigma_{\alpha},\sigma_{\gamma}\widehat{\rangle}}{\widehat{\sigma}_{\beta}
} \frac{\langle
\sigma_{\gamma},\sigma_{\alpha},\sigma_{\beta}\widehat{\rangle} }{
\widehat{\sigma}_{\gamma}} {\tau}_{\alpha,\beta;\gamma,\gamma}
\widehat{\ra}(z).
\end{equation}
\el
{\bf Proof.} There is some redundancy in  assuming that all the measures on the left hand are finite. The finiteness of some may be deduced from the rest using previous results. In (\ref{inversa2*}), the function $f$ plays the role of a wildcard  convenient in the application of that formula.

Let us prove (\ref{inversa2*}). Taking into account
(\ref{2.1}) and (\ref{2.3}), we have that
\[
\langle \langle \tau_{\alpha,\alpha}
,\sigma_{\beta},\sigma_{\alpha}\rangle,\langle  f \sigma_{\gamma},
\sigma_{\alpha}\rangle \widehat{\rangle}(z)=\langle f
\sigma_{\gamma}, \sigma_{\alpha} \widehat{\rangle}(z)\langle
\tau_{\alpha,\alpha} ,\sigma_{\beta},
{\sigma}_{\alpha}\widehat{\ra}(z) -\langle \langle f
\sigma_{\gamma}, \sigma_{\alpha}\rangle, \tau_{\alpha,\alpha}
,\sigma_{\beta},{\sigma}_{\alpha}\widehat{\ra}(z) =
\]
\[
\langle f   \sigma_{\gamma},\sigma_{\alpha}
\widehat{\rangle}(z)\left(\frac{|\langle
\sigma_{\alpha},\sigma_{\beta}\rangle|}{|\sigma_{\alpha}|}-
\frac{\langle\sigma_{\alpha},\sigma_{\beta}\widehat{\rangle}(z)}{
\widehat{\sigma}_{\alpha}(z)}\right) -
\int\left(\frac{|\langle
\sigma_{\alpha},\sigma_{\beta}\rangle|}{|\sigma_{\alpha}|}-
\frac{\langle\sigma_{\alpha},\sigma_{\beta}\widehat{\rangle}(x_{\gamma})}{
\widehat{\sigma}_{\alpha}(x_{\gamma})}\right)
\frac{f(x_{\gamma})d\langle\sigma_{\gamma}, {\sigma}_{\alpha}\ra
(x_{\gamma})}{z-x_{\gamma}}=
\]
\[
\int
\langle\sigma_{\alpha},\sigma_{\beta}\widehat{\rangle}(x_{\gamma})\frac{
f(x_{\gamma}) d\sigma_{\gamma} (x_{\gamma})}{z-x_{\gamma}}-\langle f
\sigma_{\gamma}, \sigma_{\alpha}
\widehat{\rangle}(z)\frac{\langle\sigma_{\alpha},\sigma_{\beta}\widehat{\rangle}(z)}{
\widehat{\sigma}_{\alpha}(z)}.
\]
This and (\ref{2.2}) render
\[
\frac{ \widehat{\sigma}_{\alpha}(z)}{\langle
\sigma_{\alpha},\sigma_{\beta}\widehat{\rangle}(z)} \langle \langle
\tau_{\alpha,\alpha} ,\sigma_{\beta},\sigma_{\alpha}\rangle, f
 \sigma_{\gamma}, \sigma_{\alpha} \widehat{\rangle}(z)=
\]
\[
\frac{ \widehat{\sigma}_{\alpha}(z)}{\langle
\sigma_{\alpha},\sigma_{\beta}\widehat{\rangle}(z)}\int
\langle\sigma_{\alpha},\sigma_{\beta}\widehat{\rangle}(x_{\gamma})
\frac{f(x_{\gamma})d\sigma_{\gamma}(x_{\gamma})}{z-x_{\gamma}}-\langle
f \sigma_{\gamma},\sigma_{\alpha} \widehat{\rangle}(z) =
\]
\[ -\langle
f \sigma_{\gamma},\sigma_{\alpha} \widehat{\rangle}(z)
 + \frac{|{\sigma}_{\alpha} |}{|\langle
\sigma_{\alpha},\sigma_{\beta} {\rangle}|} \la f \sigma_{\gamma},
\sigma_{\alpha}, \sigma_{\beta} \widehat{\ra} (z)  + \la f
\sigma_{\gamma}, \sigma_{\alpha}, \sigma_{\beta} \widehat{\ra} (z)
\la \frac{\tau_{\alpha,\beta}}{\widehat{\sigma}_{\beta}},\la
\sigma_{\beta},\sigma_{\alpha} \ra \widehat{\ra}(z).
\]

Fix a compact set ${\mathcal{K}} \subset {\mathbb{C}} \setminus \Delta_{\alpha}$. As in the proof of Lemma \ref{alfabeta+} take $\theta > 0$ sufficiently small so that the curve $\Gamma_{\theta}$ separates ${\mathcal{K}}$ and $\Delta_{\alpha}$. (We are considering the most complicated case when $\Delta_{\alpha}, \Delta_{\beta},$ and $\Delta_{\gamma}$ are unbounded.)

Set $w(z) = \frac{ \widehat{\sigma}_{\alpha}(z)}{\langle
\sigma_{\alpha},\sigma_{\beta}\widehat{\rangle}(z)} \langle \langle
\tau_{\alpha,\alpha} ,\sigma_{\beta},\sigma_{\alpha}\rangle, f
\sigma_{\gamma}, \sigma_{\alpha} \widehat{\rangle}(z)$.  Obviously this function satisfies c) in Lemma \ref{lem:2}. That it also fulfills d) is a consequence of the second part of Lemma \ref{lem:otro}, applied on $ \langle \langle
\tau_{\alpha,\alpha} ,\sigma_{\beta},\sigma_{\alpha}\rangle, f
\sigma_{\gamma}, \sigma_{\alpha} \widehat{\rangle}(z),$ and the first part of the last statement in Lemma \ref{alfabeta+}, employed on $\frac{ \widehat{\sigma}_{\alpha}(z)}{\langle
\sigma_{\alpha},\sigma_{\beta}\widehat{\rangle}(z)}$ (use Lebesgue's dominated convergence theorem, taking into consideration equality (\ref{2.2}), to show that the function is uniformly bounded on a sector near $x_{\alpha,\beta}$). From (\ref{eq:2.2}) it follows that
\[ w(z) = \frac{1}{2\pi i} \int_{\Gamma_{\theta}} \frac{w(\zeta)}{\zeta -z} d\zeta = \frac{1}{2\pi i} \int_{\Gamma_{\theta}} \frac{(h_1(\zeta) + g(\zeta)\widehat{s}(\zeta) )d\zeta }{\zeta -z},
\]
where $h_1(\zeta) = -\langle
f \sigma_{\gamma},\sigma_{\alpha} \widehat{\rangle}(\zeta)
 + \frac{|{\sigma}_{\alpha} |}{|\langle
\sigma_{\alpha},\sigma_{\beta} {\rangle}|} \la f \sigma_{\gamma},
\sigma_{\alpha}, \sigma_{\beta} \widehat{\ra} (\zeta) $, $g(\zeta) = \la f
\sigma_{\gamma}, \sigma_{\alpha}, \sigma_{\beta} \widehat{\ra} (\zeta)$ and $ s =
 \frac{\la
\sigma_{\beta},\sigma_{\alpha} \widehat{\ra} \tau_{\alpha,\beta}}{\widehat{\sigma}_{\beta}}$.
Using (\ref{2.1}) on $h(\zeta) = h_1(\zeta)/(\zeta -z)$ and (\ref{2.3}) on $g(\zeta)\widehat{s}(\zeta)$, (\ref{inversa2*}) readily follows. That $h, g$ and $s$ satisfy the assumptions imposed in Lemma \ref{lem:2} is a consequence of Lemma \ref{lem:otro} and the first limit in (\ref{doslim}) (notice that $s$ is the measure on the right hand of (\ref{2.2})).

In order to  prove (\ref{shirenu}), first one must derive the identity (for details see \cite[Lemma 3.3]{LF4})
\[ \frac{\langle
\sigma_{\alpha},\sigma_{\beta}\widehat{\rangle}(z)}{\langle \langle
\sigma_{\alpha}, \sigma_{\beta} \rangle,
\sigma_{\gamma}\widehat{\rangle}(z)} \la \frac{\la
\sigma_{\beta},\sigma_{\alpha}
\widehat{\ra}}{\widehat{\sigma}_{\beta}} \tau_{\alpha,\beta},
\sigma_{\gamma}, \sigma_{\alpha},\sigma_{\beta} \widehat{\ra}(z) =
\]
\[
-\frac{\langle \sigma_{\beta}, {\sigma}_{\alpha}\widehat{\ra}(z)}{
\widehat{\sigma}_{\beta}(z)}+ \frac{\la
\sigma_{\beta},\sigma_{\alpha},\sigma_{\gamma}
\widehat{\ra}(z)}{\widehat{\sigma}_{\beta}(z)}\frac{|\langle\sigma_{\alpha},\sigma_{\beta}
\rangle|}
{|\langle\langle\sigma_{\alpha},\sigma_{\beta}\rangle,\sigma_{\gamma}
{\rangle}|} + \frac{\la
\sigma_{\beta},\sigma_{\alpha},\sigma_{\gamma}
\widehat{\ra}(z)}{\widehat{\sigma}_{\beta}(z)} \la  \frac{
\tau_{\alpha,\beta;\gamma,\gamma} }{\widehat{\sigma}_{\gamma}},
\sigma_{\gamma},\sigma_{\alpha},\sigma_{\beta} \widehat{\ra}(z).
\]
Then one uses (\ref{2.1})-(\ref{2.3}), as we did above, with
\[ w(\zeta) = \frac{\langle
\sigma_{\alpha},\sigma_{\beta}\widehat{\rangle}(z)}{\langle \langle
\sigma_{\alpha}, \sigma_{\beta} \rangle,
\sigma_{\gamma}\widehat{\rangle}(\zeta)} \la \frac{\la
\sigma_{\beta},\sigma_{\alpha}
\widehat{\ra}}{\widehat{\sigma}_{\beta}} \tau_{\alpha,\beta},
\sigma_{\gamma}, \sigma_{\alpha},\sigma_{\beta} \widehat{\ra}(\zeta)
\]
\[ h_1(\zeta) = -\frac{\langle \sigma_{\beta}, {\sigma}_{\alpha}\widehat{\ra}(\zeta)}{
\widehat{\sigma}_{\beta}(\zeta)}+ \frac{\la
\sigma_{\beta},\sigma_{\alpha},\sigma_{\gamma}
\widehat{\ra}(\zeta)}{\widehat{\sigma}_{\beta}(\zeta)}\frac{|\langle\sigma_{\alpha},\sigma_{\beta}
\rangle|}
{|\langle\langle\sigma_{\alpha},\sigma_{\beta}\rangle,\sigma_{\gamma}
{\rangle}|}, \qquad h(\zeta) = \frac{h_1(\zeta)}{\zeta -z}
\]
\[ g(\zeta) = \frac{\la
\sigma_{\beta},\sigma_{\alpha},\sigma_{\gamma}
\widehat{\ra}(\zeta)}{\widehat{\sigma}_{\beta}(\zeta)}, \qquad \mbox{and} \qquad s=   \frac{\la \sigma_{\gamma},\sigma_{\alpha},\sigma_{\beta} \widehat{\ra}
 }{\widehat{\sigma}_{\gamma}}
 \tau_{\alpha,\beta;\gamma,\gamma},
\]
to conclude the proof. To justify the conditions on $w, g, h,$ and $s$ one uses arguments similar to the previous case. For example, from (\ref{2.2}) we have
\[  \frac{\la {\sigma}_{\alpha}, \sigma_{\beta} \widehat{\ra}(z)}{\la \langle
\sigma_{\alpha},\sigma_{\beta}\ra , \sigma_{\gamma} \widehat{\rangle}(z)}=
\frac{|\la {\sigma}_{\alpha}, \sigma_{\beta}\ra|}{|\la \langle
\sigma_{\alpha},\sigma_{\beta}\ra , \sigma_{\gamma}\ra|}  + \langle
 \frac{{\tau}_{\alpha,\beta;\gamma,\gamma}}{\widehat{\sigma}_{\gamma}},\sigma_{\gamma}, {\sigma}_{\alpha}, \sigma_{\beta}\widehat{\ra}(z),
\]
and due to the first limit in (\ref{doslim}) it follows that $\widehat{s}(x_{\alpha,\beta}) \in {\mathbb{R}}$. We are done.   \fp

\vspace{0,2cm}

The proof of Theorem \ref{teo:1} is based on reducing the problem to the case of multi-indices with decreasing components. We need to learn how to transpose components of the multi-index obtaining systems with the same zeros. In our aid comes the next result whose prove is based on (\ref{4.4})-(\ref{4.2}) and (\ref{inversa2*})-(\ref{shirenu}).

\bl \label{lem:4}
Let $(s_{1,1},\ldots,s_{1,m}) =
{\mathcal{N}}(\sigma_1,\ldots,\sigma_m), m \geq 1,$ be such that the first three moments of $\sigma_k, k=1,\ldots,m,$ are finite. Given ${\bf n} \in
{\mathbb{Z}}_+^{m+1}$, consider the linear form ${\mathcal{L}}_{\bf n}$ defined in Lemma $\ref{lem:3}$. Assume that
$n_j = \max\{n_0+1,n_1,\ldots,n_m\}$. Then, there exist a Nikishin system $(s_{1,1}^*,\ldots,s_{1,m}^*) =
{\mathcal{N}}(\sigma_1^*,\ldots,\sigma_m^*)$, a multi-index $\,{\bf n}^* = (n_0^*,\ldots,n_m^*) \in {\mathbb{Z}}_+^{m+1}$
which is a permutation of ${\bf n}$ with $n_0^* = n_j$, and polynomials with real coefficients $p_k^*, \deg p_k^* \leq n_k^* -1, k=0,\ldots,m$, such that
\[ {\mathcal{L}}_{\bf n} = p_0 + \sum_{k=1}^m p_k \widehat{s}_{1,k} = (p_0^* + \sum_{k=1}^m p_k^* \widehat{s}_{1,k}^*)\widehat{s}_{1,j} = {\mathcal{L}}_{\bf n}^*\widehat{s}_{1,j}.
\]
If for some $k, 1 \leq k \leq m,$ the first $n+3, n\geq 0,$ moments of $\sigma_k$ exist, then the first $n+1$ moments of $\sigma_k^*$ also exist.
\el
\vspace{0,2cm}

{Proof.} Assume that $j=1$. From (\ref{eq:1}) and (\ref{4.4}), we have
\[ \frac{ {\mathcal{L}}_{\bf n}}{\widehat{s}_{1,1} } =
\frac{ p_0}{ \widehat{s}_{1,1}}+p_1   + \sum_{k=2}^m p_k\frac{\widehat{s}_{1,k}}{\widehat{s}_{1,1}} =
\]
\[
 (\ell_{1,1}p_0  + p_1 + \sum_{k=2}^m \frac{|s_{1,k}|}{|s_{1,1}|}p_k) +  p_0\widehat{\tau}_{1,1}
- \sum_{k=2}^m p_k \la  {\tau}_{1,1},
s_{2,k},\sigma_1 \widehat{\ra} = {\mathcal{L}}_{\bf n}^*.
\]
We are done taking ${\bf n}^* = (n_1,n_0,n_2,\ldots,n_m)$ and
\[{\mathcal{N}}(\sigma_1^*,\ldots,\sigma_m^*) = {\mathcal{N}}(\tau_{1,1},\la \sigma_2,\sigma_1\ra,\sigma_3,\ldots,\sigma_m)
\]
since $\la s_{2,k},\sigma_1\ra = \la \la \sigma_2,\sigma_1\ra ,\sigma_3,\ldots,\sigma_k\ra$ when $k\geq 3$.  \hfill $\Box$ \vspace{0,2cm}

In the sequel $2 \leq j \leq m$. From (\ref{eq:1}),  (\ref{4.5}), and the first equalities in
(\ref{4.6})-(\ref{4.2}), one has
\[ \frac{ {\mathcal{L}}_{\bf n}}{\widehat{s}_{1,j}} =
\frac{p_0}{\widehat{s}_{1,j}}+
 p_j   + \sum_{k\neq j,k=1}^m p_k  \frac{\widehat{s}_{1,k}}{\widehat{s}_{1,j}} =
 (\ell_{1,j}p_0  +  p_j  + \sum_{k\neq j, k=1}^m \frac{|s_{1,k}|}{|s_{1,j}|}
p_k)  +
\]
\[p_0 \widehat{\tau}_{1,j}
+ p_1\la
\frac{\tau_{1,j}}{\widehat{s}_{2,j}}, \la s_{2,j},\sigma_1 \ra
\widehat{\ra} +
\]
\[ \sum_{k=2}^{j-1} (-1)^{k-1} p_k \la
\tau_{1,j},\la\tau_{2,j},s_{1,j} \ra,\ldots,
\la\tau_{k-1,j},s_{k-2,j}\ra, \frac{\la \tau_{k,j},
{s}_{k-1,j}\ra}{\widehat{s}_{k+1,j}},\la s_{k+1,j},\sigma_k\ra
\widehat{\ra} +
\]
\begin{equation}
\label{eq:L}
(-1)^j  \sum_{k=j+1}^{m} p_k \la
\tau_{1,j},\la\tau_{2,j},s_{1,j} \ra,\ldots,
\la\tau_{j,j},s_{j-1,j}\ra,\la s_{j+1,k},\sigma_j\ra \widehat{\ra}.
\end{equation}
Now, it is not so clear who the auxiliary Nikishin system should be because some annoying ratios of Cauchy transforms have appeared. We shall see that already for $j=2$ there are two candidates, and for general $j$ the number of candidates equals $2^{j-1}$.

We can use (\ref{4.8}) (see the second inequalities in (\ref{4.6})-(\ref{4.2})) to obtain
\[  \frac{ {\mathcal{L}}_{\bf n}}{\widehat{s}_{1,j}}  = (\ell_{1,j}p_0  +  p_j  + \sum_{k\neq j, k=1}^m \frac{|s_{1,k}|}{|s_{1,j}|}
p_k) + (p_0 + \frac{|\la
 {s}_{2,j},\sigma_1  \ra |}{|{s}_{2,j}|}p_1)\tau_{1,j} +
\]
\[ \sum_{k=2}^{j-1}
 (-1)^{k-1} (p_{k-1} + \frac{|\la
 {s}_{k+1,j},\sigma_k  \ra |}{|{s}_{k+1,j}|} p_k) \la \tau_{1,j},\la
 \tau_{2,j},s_{1,j}\ra,\ldots,\la \tau_{k,j},s_{k-1,j} \ra
 \widehat{\ra} +
\]
\[ (-1)^{j-1} p_{j-1} \la \tau_{1,j}, \la
 \tau_{2,j},s_{1,j}\ra,\ldots,\la \tau_{j,j},s_{j-1,j} \ra
 \widehat{\ra} +
\]
\begin{equation} \label{4.9}
(-1)^j \sum_{k=j+1}^{m} p_k \la \tau_{1,j},\la\tau_{2,j},s_{1,j} \ra,\ldots,
\la\tau_{j,j},s_{j-1,j}\ra,\la s_{j+1,k},\sigma_j\ra \widehat{\ra}.
\end{equation}
(The sum $\sum_{k=2}^{j-1}$ is empty if $j=2$.)

If we are in the class  of multi-indices
\[ {\mathbb{Z}}_+^{m+1}(*) = \{{\bf n} \in {\mathbb{Z}}_+^{m+1}: \not\exists \,\, 0 \leq i < j < k < m\,\, \mbox{such that}\,\, n_i <  n_j \leq n_k \},
\]
and take $j$ to be the first component for which $n_j = \max\{n_0+1,n_1,\ldots,n_m\}$, then  $n_0 \geq \cdots \geq n_{j-1}$. It follows that
\[ \deg (\ell_{1,j}p_0  +  p_j  + \sum_{k\neq j, k=1}^m \frac{|s_{1,k}|}{|s_{1,j}|}
p_k) \leq n_j -1
\]
and
\[ \deg (p_{k-1} + \frac{|\la
 {s}_{k+1,j},\sigma_k  \ra |}{|{s}_{k+1,j}|} p_k) \leq n_{k-1} -1,\qquad k=1,\ldots,j-1.
\]
Thus ${\mathcal{L}}_{\bf n}^*$ is the right hand side of (\ref{4.9}), which is a linear form generated by the multi-index ${\bf n}^* = (n_j,n_0\ldots,n_{j-1},n_{j+1},\ldots,n_m) \in {\mathbb{Z}}_+^{m+1}$ and the Nikishin system
\[{\mathcal{N}}(\sigma_1^*,\ldots,\sigma_m^*) = {\mathcal{N}}(\tau_{1,j},\la\tau_{2,j},s_{1,j} \ra,\ldots,
\la\tau_{j,j},s_{j-1,j}\ra,\la \sigma_{j+1},\sigma_j\ra, \sigma_{j+2},\ldots,\sigma_m).
\]

This would be sufficient to prove the AT property within the class ${\mathbb{Z}}_+^{m+1}(*)$ because   then $(n_0,\ldots,n_{j-1},n_{j+1},\ldots,n_m) \in {\mathbb{Z}}_+^m(*)$ and one can use induction. When the supports of the generating measures are non intersecting and bounded, this result was given in \cite[Theorem 2]{LIF}.

Of course, (\ref{4.9}) is still valid in the general case but, if it is not true that $n_0 \geq \ldots \geq n_{j-1}$, some of the degrees of the polynomials in the linear form on the right hand blow up with respect to the bounds established by the components of ${\bf n}^*$. We must proceed with caution.

Set
\[ {\mathcal{L}}_{\bf n}^* = p_0^* + p_0 \widehat{\tau}_{1,j}
+ p_1\la
\frac{\tau_{1,j}}{\widehat{s}_{2,j}}, \la s_{2,j},\sigma_1 \ra
\widehat{\ra}  +
\]
\[ \sum_{k=2}^{j-1} (-1)^{k-1} p_k \la
\tau_{1,j},\la\tau_{2,j},s_{1,j} \ra,\ldots,
\la\tau_{k-1,j},s_{k-2,j}\ra, \frac{\la \tau_{k,j},
{s}_{k-1,j}\ra}{\widehat{s}_{k+1,j}},\la s_{k+1,j},\sigma_k\ra
\widehat{\ra} +
\]
\begin{equation}
\label{eq:A}
(-1)^j  \sum_{k=j+1}^{m} p_k \la
\tau_{1,j},\la\tau_{2,j},s_{1,j} \ra,\ldots,
\la\tau_{j,j},s_{j-1,j}\ra,\la s_{j+1,k},\sigma_j\ra \widehat{\ra},
\end{equation}
where
\[ p_0^* = \ell_{1,j}p_0  +  p_j  + \sum_{k\neq j, k=1}^m \frac{|s_{1,k}|}{|s_{1,j}|}
p_k, \quad \deg p_0^* \leq n_j -1 = n_0^* -1.
\]
We took this function from the right hand side of (\ref{eq:L}). We must show that there exist a multi-index ${\bf n}^* \in {\mathbb{Z}}_+^{m+1}$, which is a permutation of ${\bf n}$,
and a Nikishin system ${\mathcal{N}}(\sigma_0^*,\ldots,\sigma_m^*)$ which allow to express ${\mathcal{L}}_{\bf n}^*$ as a linear form generated by them with polynomials with real coefficients. So far, $n_0^*$ defined above serves the purpose of being the first component of ${\bf n}^*$ and $p_0^*$ of being the polynomial part of the linear form.

{\bf First step.} We ask whether  $n_0 \geq n_1$ or $n_0 \leq n_1$? (When $n_0 =
n_1$ we can proceed either ways.)

\noindent {\bf A1)} If $n_0 \geq n_1$, take ${n}_1^* = n_0$  and $\sigma_1^* = \tau_{1,j}$.
Decompose $\frac{ \la s_{2,j},\sigma_1
\widehat{\ra}}{\widehat{s}_{2,j}}$ using (\ref{2.3}). Then, the
first three terms of ${\mathcal{L}}_{\bf n}^*$ are
\[ p_0^* + p_0 \widehat{\sigma}_{1}^*
+ p_1\la
\frac{\sigma_{1}^*}{\widehat{s}_{2,j}}, \la s_{2,j},\sigma_1 \ra
\widehat{\ra}  = p_0^* + (p_0
+ \frac{|\la
 {s}_{2,j},\sigma_1  \ra |}{|{s}_{2,j}|}p_1) \widehat{\sigma}_{1}^* - p_1 \la {\sigma}_{1}^*, \la \tau_{2,j},
s_{1,j} \ra \widehat{\ra}.
\]
Consequently, taking  $ {p}_1^* =  p_0 + \frac{|\la
 {s}_{2,j},\sigma_1  \ra |}{|{s}_{2,j}|}p_1$, we have that $\deg  {p}_1^* \leq  {n}_1^* -1 (= n_0-1)$.

In case that $j=2$, we obtain
\[ {\mathcal{L}}_{\bf n}^* = p_0^* +  {p}_1^* \widehat{\tau}_{1,2} - p_1 \la \tau_{1,2},
 \la \tau_{2,2},s_{1,2} \ra \widehat{\ra} + \sum_{k=3}^{m} p_k
\la \tau_{1,2}, \la\tau_{2,2},s_{1,2} \ra, \la s_{3,k},\sigma_2 \ra
\widehat{\ra}
\]
(compare with (\ref{4.9})). Then, the proof would be complete taking
$ {\bf n}^* = (n_2, n_0,n_1,n_3\ldots,n_m)$ and the Nikishin
system
\[{\mathcal{N}}(\sigma_1^*,\ldots,\sigma_m^*) = {\mathcal{N}}(\tau_{1,2},\la \tau_{2,2},s_{1,2}\ra, \la
\sigma_3,\sigma_2 \ra, \sigma_4,\ldots,\sigma_m ).
\]
(If $m=2$, then
$ {\bf n}^* = (n_2, n_0,n_1)$ and the Nikishin system is
${\mathcal{N}}(\tau_{1,2}, \la\tau_{2,2},s_{1,2}\ra)$.)

If $j \geq 3$, we obtain
\[ {\mathcal{L}}_{\bf n}^* = p_0^* + {p}_1^* \widehat{\sigma}_{1}^* - p_1 \la
 \sigma_{1}^*, \la\tau_{2,j},s_{1,j} \ra \widehat{\ra} - p_2
\la \sigma_{1}^*, \frac{\la s_{3,j},\sigma_2  \widehat{\ra}}{ \widehat{s}_{3,j} }
\la \tau_{2,j}, s_{1,j} \ra \widehat{\ra} +
\]
\[\sum_{k=3}^{j-1} (-1)^{k-1} p_k \la \sigma_{1}^*, \la\tau_{2,j},s_{1,j}
\ra,\ldots, \la\tau_{k-1,j},s_{k-2,j}\ra, \frac{\la \tau_{k,j},
{s}_{k-1,j}\ra}{\widehat{s}_{k+1,j}},\la s_{k+1,j},\sigma_k\ra
\widehat{\ra}+   \]
\[  (-1)^j \sum_{k=j+1}^{m} p_k
\la \sigma_{1}^*,\la\tau_{2,j},s_{1,j} \ra,\ldots, \la\tau_{j,j},s_{j-1,j}\ra,\la
s_{j+1,k},\sigma_j\ra \widehat{\ra}.
\]

\noindent {\bf B1)} If $n_0 \leq n_1$, take $ {n}_1^* = n_1$ and $\sigma_1^* =  \frac{\la s_{2,j},\sigma_1\widehat{\ra}}{\widehat{s}_{2,j}} \tau_{1,j}$. We can rewrite (\ref{eq:A}) as follows
\[ {\mathcal{L}}_{\bf n}^* = p_0^* + p_1\widehat{\sigma}_1^*  + p_0 \la \frac{ \widehat{s}_{2,j}}{\la s_{2,j},\sigma_1
\widehat{\ra}} {\sigma}_1^*\widehat{\ra}
+
\]
\[ \sum_{k=2}^{j-1} (-1)^{k-1} p_k \la
\frac{ \widehat{s}_{2,j}}{\la s_{2,j},\sigma_1
\widehat{\ra}} {\sigma}_1^*,\la\tau_{2,j},s_{1,j} \ra,\ldots,
\la\tau_{k-1,j},s_{k-2,j}\ra, \frac{\la \tau_{k,j},
{s}_{k-1,j}\ra}{\widehat{s}_{k+1,j}},\la s_{k+1,j},\sigma_k\ra
\widehat{\ra} +
\]
\begin{equation}
\label{eq:B}
(-1)^j  \sum_{k=j+1}^{m} p_k \la
\frac{ \widehat{s}_{2,j}}{\la s_{2,j},\sigma_1
\widehat{\ra}} {\sigma}_1^*,\la\tau_{2,j},s_{1,j} \ra,\ldots,
\la\tau_{j,j},s_{j-1,j}\ra,\la s_{j+1,k},\sigma_j\ra \widehat{\ra},
\end{equation}

Decompose $\frac{ \widehat{s}_{2,j}}{\la s_{2,j},\sigma_1
\widehat{\ra}}$ using (\ref{2.2}). Then, the
first three terms of ${\mathcal{L}}_{\bf n}^*$ in (\ref{eq:B}) can be expressed as
\[ p_0^* + p_1 \widehat{\sigma}_{1}^*
+ p_0 \la \frac{ \widehat{s}_{2,j}}{\la s_{2,j},\sigma_1
 \widehat{\ra}}  {\sigma}_{1}^* \widehat{\ra} = p_0^* + (p_1
+ \frac{|{s}_{2,j}|}{|\la
 {s}_{2,j},\sigma_1  \ra |}p_0) \widehat{\sigma}_{1}^* + p_0 \la \sigma_1^*, \frac{\widehat{s}_{1,j}}{\widehat{\sigma}_1}  \tau_{2,j;1,1} \widehat{\ra}.
\]
Taking  $ {p}_1^* = p_1 + \frac{|{s}_{2,j}|}{|\la
 {s}_{2,j},\sigma_1  \ra |}p_0  $, we have that $\deg  {p}_1^* \leq  {n}_1^* -1 (= n_1-1)$.

If $j=2$, due to (\ref{inversa2*}) (in the next formula
$\widehat{s}_{4,k} \equiv 1$ if $k=3$)
\[ \frac{\widehat{s}_{2,2}}{\la s_{2,2},\sigma_1 \widehat{\ra}}
\la \la\tau_{2,2},s_{1,2} \ra,\la s_{3,k},\sigma_2\ra \widehat{\ra}
= \frac{\widehat{\sigma}_{2}}{\la \sigma_{2},\sigma_1 \widehat{\ra}}
\la \la\tau_{2,2},s_{1,2} \ra,\widehat{s}_{4,k} \sigma_{3},\sigma_2
\widehat{\ra} =
\]
\[ \la \frac{\la \sigma_{1},\sigma_{2}
\widehat{\ra}}{\widehat{\sigma}_{1}} \tau_{2,2;1,1},
\widehat{s}_{4,k} \sigma_{3}, \sigma_{2},\sigma_{1}  \widehat{\ra} =
\la \frac{\la \sigma_{1},\sigma_{2}
\widehat{\ra}}{\widehat{\sigma}_{1}} \tau_{2,2;1,1}, \la \sigma_{3},
\sigma_{2},\sigma_{1} \ra, s_{4,k}\widehat{\ra}.
\]
Consequently,
\[ {\mathcal{L}}_{\bf n}^* = p_0^* + p_1^* \widehat{\sigma}_{1}^* + p_0 \la \sigma_1^*, \frac{\widehat{s}_{1,2}}{\widehat{\sigma}_1}  \tau_{2,2;1,1} \widehat{\ra} + \sum_{k=3}^m p_k \la \sigma_1^*, \frac{\widehat{s}_{1,2}}{\widehat{\sigma}_1}  \tau_{2,2;1,1} , \la \sigma_{3},
\sigma_{2},\sigma_{1} \ra, s_{4,k}\widehat{\ra}.
\]
In this situation, we would be done considering $ {\bf n}^* =
(n_2, n_1,n_0,n_3,\ldots,n_m)$ and the system
\[{\mathcal{N}}(\sigma_1^*,\ldots,\sigma_m^*) = {\mathcal{N}}(\frac{\la \sigma_{2},\sigma_1\widehat{\ra}}{\widehat{\sigma}_{2}} \tau_{1,2},\frac{\la \sigma_{1},\sigma_{2}
\widehat{\ra}}{\widehat{\sigma}_{1}} \tau_{2,2;1,1}, \la \sigma_{3},
\sigma_{2},\sigma_{1} \ra, \sigma_{4},\ldots,\sigma_m).
\]
(Should
$m=2$, then $ {\bf n}^* = (n_2,n_1,n_0)$ and the Nikishin system
is ${\mathcal{N}}(\frac{\la \sigma_{2},\sigma_1\widehat{\ra}}{\widehat{\sigma}_{2}} \tau_{1,j},\frac{\la \sigma_{1},\sigma_{2}
\widehat{\ra}}{\widehat{\sigma}_{1}} \tau_{2,2;1,1})$.) Therefore, if $j=2$ we are done.

Let us assume that $j \geq 3$. (The algorithm ends after $j-1$
steps.)  So far, we have used the notation $s_{k,l}$ only with $k \leq l$. We will extend its meaning to $k > l$ in which case
\[s_{k,l} = \la \sigma_k, \sigma_{k-1}, \ldots, \sigma_l \ra, \qquad k > l.
\]
Notice that if $l < k < j$, then
\[ \la s_{k,j}, s_{k-1, l} \ra = \la s_{k,l}, s_{k+1,j} \ra.
\]
The inverse measure of $\la s_{k,j}, s_{k-1, l} \ra$ we denote by $\tau_{k,j;k-1,l}$; that is,
\[ {1}/{\la s_{k,j}, s_{k-1, l} \widehat{\ra}} = \ell_{k,j;k-1,l} +  \widehat{\tau}_{k,j;k-1,l}
\]
In particular, $\tau_{2,j;1,1}$ denotes the inverse measure of $\la s_{2,j},\sigma_1 \ra$.

Let us transform the measures in $\sum_{k=2}^{j-1}$ of (\ref{eq:B}). Regarding the term with $p_2$, using (\ref{shirenu}) with
$\sigma_{\alpha} = \sigma_2, \sigma_{\beta} = s_{3,j}$ and
$\sigma_{\gamma} = \sigma_1$, we obtain
\[ \frac{\widehat{s}_{2,j}}{\la s_{2,j},\sigma_1
\widehat{\ra}}\la \frac{\la \tau_{2,j}, s_{1,j}\ra}{
\widehat{s}_{3,j} }, \la s_{3,j},\sigma_2 \ra \widehat{\ra} = \la
\frac{\la s_{3,j},\sigma_2,\sigma_1\widehat{\ra}}{\widehat{s}_{3,j}}
\frac{\widehat{s}_{1,j}}{\widehat{\sigma}_1}\tau_{2,j;1,1}
\widehat{\ra}.
\]
For $j=3$, $\sum_{k=3}^{j-1}$ is empty, so here the formulas make
sense when $j\geq 4$. Using (\ref{inversa2*}), with $\sigma_{\alpha}
= s_{2,j},$ $\sigma_{\beta} = \sigma_1,$ $\sigma_{\gamma} =
\tau_{3,j}$, and
\[ f = f_{1,j,k} = \left\{
\begin{array}{ll}
\displaystyle{\frac{\la s_{4,j}, \sigma_3 \widehat{\ra}}{\widehat{s}_{4,j}}}, &  3 = k < j \leq m, \\
\la \displaystyle{\frac{\la
\tau_{4,j},s_{3,j}\ra}{\widehat{s}_{5,j}}},\la s_{5,j},\sigma_4 \ra
\widehat{\ra}, &  4 = k < j \leq m, \\
\la \la \tau_{4,j},s_{3,j}\ra,\ldots,\la \tau_{k-1,j},s_{k-2,j} \ra,
\displaystyle{\frac{\la \tau_{k,j},s_{k-1,j}
\ra}{\widehat{s}_{k+1,j}}} , \la s_{k+1,j},\sigma_k \ra
\widehat{\ra}, &  5 \leq k < j \leq m,
\end{array}
\right.
\]
we obtain
\[ \frac{\widehat{s}_{2,j}}{\la s_{2,j},\sigma_1 \widehat{\ra}} \la \la\tau_{2,j},s_{1,j}
\ra,\ldots, \la\tau_{k-1,j},s_{k-2,j}\ra, \frac{\la \tau_{k,j},
{s}_{k-1,j}\ra}{\widehat{s}_{k+1,j}},\la s_{k+1,j},\sigma_k\ra
\widehat{\ra} =
\]
\[ \la  {\frac{\widehat{s}_{1,j}}{\widehat{\sigma}_1}} \tau_{2,j;1,1} ,
f_{1,j,k} \tau_{3,j},s_{2,j},\sigma_1  \widehat{\ra} =
\]
\[
\left\{
\begin{array}{ll}
\la \displaystyle{\frac{\widehat{s}_{1,j}}{\widehat{\sigma}_1}}
\tau_{2,j;1,1} , \displaystyle{\frac{\la \tau_{3,j},s_{2,j},\sigma_1
\ra }{\widehat{s}_{4,j}}}, \la s_{4,j},\sigma_3 \ra
\widehat{\ra}, &  3 = k < j \leq m, \\
\la \displaystyle{\frac{\widehat{s}_{1,j}}{\widehat{\sigma}_1}}
\tau_{2,j;1,1} , \la \tau_{3,j},s_{2,j},\sigma_1 \ra,\displaystyle{
\frac{\la \tau_{4,j}, s_{3,j} \ra}{\widehat{s}_{5,j}} },\la s_{5,j},
\sigma_4 \ra \widehat{\ra}, &  4 = k < j \leq m,\\
\la  {\frac{\widehat{s}_{1,j}}{\widehat{\sigma}_1}} \tau_{2,j;1,1} ,
\la \widehat{\sigma}_1 s_{2,j} \widehat{\ra}\tau_{3,j}, \widehat{s}_{3,j} \tau_{4,j},
\ldots, \widehat{s}_{k-2,j}\tau_{k-1,j} ,  {\frac{\la
\widehat{\sigma}_k s_{k+1,j} \widehat{\ra}\widehat{s}_{k-1,j} }{\widehat{s}_{k+1,j}}}
\tau_{k,j}  \widehat{\ra}, & 5 \leq k < j \leq m.
\end{array}
\right.
\]
In the last row, we used a more compact notation to fit the line. Notice that it is the same as
\[ \la  {\frac{\widehat{s}_{1,j}}{\widehat{\sigma}_1}} \tau_{2,j;1,1} ,
\la \tau_{3,j}, s_{2,j},\sigma_1  {\ra}, \la \tau_{4,j},  {s}_{3,j} \ra,
\ldots, \la \tau_{k-1,j},s_{k-2,j}\ra,  {\frac{\la
\tau_{k,j},s_{k-1,j} \ra}{\widehat{s}_{k+1,j}}}, \la
s_{k+1,j},\sigma_k \ra \widehat{\ra}.
\]

As for the terms in  $\sum_{k = j+1}^m$, applying
(\ref{inversa2*}) with $\sigma_{\alpha} = s_{2,j},$ $\sigma_{\beta}
= \sigma_1,$ $\sigma_{\gamma} = \tau_{3,j}$, and
\[ f  = f_{1,j,k}= \left\{
\begin{array}{ll}
\la s_{4,k}, \sigma_3 \widehat{\ra},  & 3 = j < k \leq m, \\
\la \la \tau_{4,j},s_{3,j}\ra,\ldots,\la \tau_{j,j},s_{j-1,j}
\ra,\la s_{j+1,k},\sigma_j \ra  \widehat{\ra}, &  4 \leq j < k \leq
m,
\end{array}
\right.
\]
it follows that
\[ \frac{\widehat{s}_{2,j}}{\la s_{2,j},\sigma_1 \widehat{\ra}}
\la \la\tau_{2,j},s_{1,j} \ra,\ldots, \la\tau_{j,j},s_{j-1,j}\ra,\la
s_{j+1,k},\sigma_j\ra \widehat{\ra} =
\]
\[ = \la \displaystyle{\frac{\widehat{s}_{1,j}}{\widehat{\sigma}_1}}
\tau_{2,j;1,1} , f_{1,j,k} \tau_{3,j},s_{2,j},\sigma_1 \widehat{\ra}
=
\]
\[
\left\{
\begin{array}{ll}
\la \displaystyle{\frac{\widehat{s}_{1,3}}{\widehat{\sigma}_1}}
\tau_{2,3;1,1} , \la \tau_{3,3},s_{2,3},\sigma_1 \ra, \la
s_{4,k},\sigma_3 \ra
\widehat{\ra},  & 3 = j < k \leq m, \\
\la \displaystyle{\frac{\widehat{s}_{1,j}}{\widehat{\sigma}_1}}
\tau_{2,j;1,1} , \la \tau_{3,j},s_{2,j},\sigma_1 \ra,\la \tau_{4,j},
s_{3,j} \ra, \ldots, \la \tau_{j,j},s_{j-1,j}\ra, \la s_{j+1,k},
\sigma_j \ra \widehat{\ra},  & 4 \leq j < k \leq m.
\end{array}
\right.
\]
When $j=m$ no such terms exist.

Therefore,
\[
 {\mathcal{L}}_{\bf n}^* = p_0^* + p_1^* \widehat{\sigma}_{1}^* + p_0 \la \sigma_1^*, \frac{\widehat{s}_{1,j}}{\widehat{\sigma}_1}  \tau_{2,j;1,1} \widehat{\ra}
  - p_2 \la \sigma_1^*, \frac{ \la
s_{3,j},\sigma_2,\sigma_1\widehat{\ra}}{\widehat{s}_{3,j}}
\frac{\widehat{s}_{1,j}}{\widehat{\sigma}_1}\tau_{2,j;1,1}\widehat{\ra}
+
\]
\[ \sum_{k=3}^{j-1} (-1)^{k-1} p_k \la \sigma_1^*,\displaystyle{\frac{\widehat{s}_{1,j}}{\widehat{\sigma}_1}}
\tau_{2,j;1,1} , f_{1,j,k} \tau_{3,j},s_{2,j},\sigma_1 \widehat{\ra}
+ (-1)^j \sum_{k= j+1}^m p_k \la \sigma_1^*,
\displaystyle{\frac{\widehat{s}_{1,j}}{\widehat{\sigma}_1}}
\tau_{2,j;1,1} , f_{1,j,k} \tau_{3,j},s_{2,j},\sigma_1
\widehat{\ra}.
\]
Using the notation for $f_{1,j,k}$ defined previously, in {\bf A1)} we ended up with
\[ {\mathcal{L}}_{\bf n}^* = p_0^* + {p}_1^* \widehat{\sigma}_{1}^* - p_1 \la
 \sigma_{1}^*, \la\tau_{2,j},s_{1,j} \ra \widehat{\ra} - p_2
\la \sigma_{1}^*, \frac{\la s_{3,j},\sigma_2  \widehat{\ra}}{ \widehat{s}_{3,j} }
\la \tau_{2,j}, s_{1,j} \ra \widehat{\ra} +
\]
\[\sum_{k=3}^{j-1} (-1)^{k-1} p_k \la \sigma_{1}^*, \la\tau_{2,j},s_{1,j}
\ra,f_{1,j,k} \tau_{3,j}, s_{2,j} \widehat{\ra} +
  (-1)^j \sum_{k=j+1}^{m} p_k
\la \sigma_{1}^*,\la\tau_{2,j},s_{1,j} \ra,f_{1,j,k} \tau_{3,j}, s_{2,j} \widehat{\ra}.
\]
Denote
\[  \sigma_2^{(l_1)} = \left\{
\begin{array}{ll}
\widehat{s}_{1,j}\tau_{2,j},  & l_1=1, \\
\frac{\widehat{s}_{1,j}}{\widehat{\sigma}_1}\tau_{2,j;1,1}, & l_1=0.
\end{array}
\right.
\]
The two formulas for ${\mathcal{L}}_{\bf n}^*$ may be put in one writing
\begin{equation}
\label{eq:C}
{\mathcal{L}}_{\bf n}^* = p_0^* + {p}_1^* \widehat{\sigma}_{1}^* + (-1)^{l_1} p_{l_1} \la
 \sigma_{1}^*, \sigma_2^{(l_1)} \widehat{\ra} - p_2
\la \sigma_{1}^*, \frac{\la s_{3,j},s_{2,l_1+1}  \widehat{\ra}}{ \widehat{s}_{3,j} }
\sigma_2^{(l_1)}\widehat{\ra} +
\end{equation}
\[\sum_{k=3, k\neq j}^{m} \delta_{k,j} p_k \la \sigma_{1}^*, \sigma_2^{(l_1)},f_{1,j,k} \tau_{3,j}, s_{2,l_1+1},s_{3,j} \widehat{\ra},  \qquad n_{l_1} = \min\{n_0,n_1\},
\]
where
\[ \delta_{k,j} = \left\{
\begin{array}{cc}
(-1)^{k-1}, & k < j, \\
(-1)^{j}, & k > j.
\end{array}
\right.
\]
We also have
\[ \deg p_0^* \leq n_j -1 = n_0^* -1, \qquad \deg p_1^* \leq \max\{n_0,n_1\} -1 = n_1^* -1.
\]

For $j=2$ we found the following solutions
\[
\begin{array}{lll}
\underline{\max\{n_0,n_1\}} & \underline{(n_0^*,\ldots,n_m^*)} & \underline{{\mathcal{N}}(\sigma_1^*,\ldots,\sigma_m^*)} \\
n_0 & (n_2,n_0,n_1,n_3,\ldots,n_m) & {\mathcal{N}}(\tau_{1,2},\la \tau_{2,2},s_{1,2}\ra, s_{3,2}, \sigma_4,\ldots,\sigma_m ) \\
n_1 & (n_2,n_1,n_0,n_3,\ldots,n_m) & {\mathcal{N}}(\frac{\widehat{s}_{2,1}}{\widehat{s}_{2,2}} \tau_{1,2},\frac{\widehat{s}_{1,2}}{\widehat{s}_{1,1}} \tau_{2,2;1,1}, s_{3,1}, \sigma_{4},\ldots,\sigma_m)
\end{array}
\]
(If $m=2$, then $ {\bf n}^*$ has only the first three components and
the Nikishin system has only the first two measures indicated.)

The case $j=2$ has been solved; therefore $j\geq 3$. When $j=3$, one can further reduce (\ref{eq:C}) following the arguments employed in {\bf A1)} should $\min\{n_0,n_1\} \geq n_2$ or in {\bf B1)} in case that
$\min\{n_0,n_1\} \leq n_2$. (When $\min\{n_0,n_1\} = n_2$ one can
proceed either ways.) To complete the proof of the formula one uses induction on $j$. For details see \cite[Lemma 2.3]{LF4}. In particular, there you will find tables with all 4 solutions for $j=3$ and the eight solutions for $j=4$.

Taking into consideration the general expression of the measures $\sigma_k^*, k=1,\ldots,m,$ (for this see the induction proof of \cite[Lemma 2.3]{LF4}), if the first $n+3, n \geq 0,$ moments of $\sigma_k$ exist, from Lemma \ref{momentos} (see also Remark \ref{rem:1}) and the second limit in (\ref{doslim}), it follows that the first $n+1$ moments of $\sigma_k^*$ are finite. With this we conclude the proof. \hfill $\Box$
\vspace{0,2cm}

\section{Proof of main results and other consequences.}

{\bf Proof of Theorem \ref{teo:1}.} Obviously, the first statement of the theorem follows from the second. We prove the second one using induction on $m$. For $m=0$ the linear form reduces to a polynomial of degree $\leq n_0-1$ and thus has at most $n_0 -1$ zeros in the complex plane as claimed.

Assume that the result is true for any Nikishin system with $m-1 (\geq 0) $ measures satisfying the assumptions of Theorem \ref{teo:1} and let us show that it is also valid for Nikishin systems with $m$ measures. To the contrary, let us suppose that ${\mathcal{L}}_{\bf n}$ has at least $|{\bf n}|$ zeros on ${\mathbb{C}}\setminus \Delta_1$.

Should $n_0 = \max\{n_0,n_1,\ldots,n_m\}$, by Lemma \ref{lem:3} the linear form $p_1 + \sum_{k=2}^m p_k \widehat{s}_{2,k}$ would have at least $|{\bf n}| - n_0$ zeros in ${\mathbb{C}} \setminus \Delta_2$. Now, $|{\bf n}| - n_0$ is the norm of the multi-index $(n_1,\ldots,n_m)$ which together with the Nikishin system ${\mathcal{N}}(\sigma_2,\ldots,\sigma_m)$ define the reduced form. This contradicts the induction hypothesis.

Suppose that $n_j = \max\{n_0+1,n_1,\ldots,n_m\}$. According to Lemma \ref{lem:4}, the linear form ${\mathcal{L}}_{\bf n}^* $ has the same zeros as ${\mathcal{L}}_{\bf n}$ in ${\mathbb{C}}\setminus \Delta_1$, since $\widehat{s}_{1,j}$ is never zero on that region. The multi-index ${\bf n}^*$ which determines ${\mathcal{L}}_{\bf n}^*$ has the same norm as ${\bf n}$  and its first component satisfies the assumptions of Lemma $\ref{lem:3}$. Thus, the reduced form
$(p_1^* + \sum_{k=2}^m p_k^* \widehat{s}_ {2,k}^*)$ has at least $|{\bf n}| - n_j = |{\bf n}^*| - n_0^*$ zeros in $\Delta_1 \subset {\mathbb{C}}\setminus \Delta_2$ contradicting the induction hypothesis. In fact, ${\mathcal{L}}_{\bf n}^*$ is generated by the multi-index $(n_1^*,\ldots,n_m^*)$ and ${\mathcal{N}}(\sigma_2^*,\ldots,\sigma_m^*)$ whose $m-1$ generating measures have the required number of finite moments according to the last statement of Lemma \ref{lem:4}. We are done. \hfill $\Box$ \vspace{0,2cm}

{\bf Proof of Theorem \ref{lem:orto}.}. Obviously, by linearity,
(\ref{orto}) is equivalent to showing that for each $j=0,\ldots,m_2$ and each $\nu =0,\ldots,n_{2,j} -1$
\[ \int x^{\nu} \widehat{s}_{1,j}^2(x) \sum_{k=0}^{m_2} a_{{\bf n},k}(x) \widehat{s}_{1,k}^1(x) d\sigma_0^2(x) =0,
\]
where $\widehat{s}_{1,0}^2 \equiv 1.$ According to b) of Definition \ref{defiSor}, we know that
\[ z^{\nu}\left( \sum_{k=0}^{m_2} a_{{\bf n},k}(z) \int \frac{\widehat{s}_{1,j}^2(x) \widehat{s}_{1,k}^1(x)d\sigma_0^2(x)}{z-x} - d_{{\bf n},j}(z)\right) = {\mathcal{O}}\left(\frac{1}{z^2}\right), \qquad z \to \infty.
\]
In particular, this implies that the coefficient corresponding to $1/z^{\nu +1}$ in the asymptotic expansion at $\infty$ of the function in parenthesis on the left hand side equals zero.

For each fixed $j$ and $k$ denote
\[ c_{j,k,n} = \int x^{n} \widehat{s}_{1,j}^2(x) \widehat{s}_{1,k}^1(x)d\sigma_0^2(x)
\]
and let
\[ a_{{\bf n},k}(z) = a_{{\bf n},k,0} + a_{{\bf n},k,1}z + \cdots + a_{{\bf n},k,n_{1,k}-1} z^{n_{1,k}-1}.
\]
Taking the asymptotic expansion at $\infty$ of the integrals inside the parenthesis  above, it is easy to verify that the coefficient of $1/z^{\nu +1}$ equals
\[ \sum_{k=0}^{m_2}  \sum_{r=0}^{n_{1,k} -1} a_{{\bf n},k,r} c_{j,k,\nu +r} = \sum_{k=0}^{m_2}  \sum_{r=0}^{n_{1,k} -1} a_{{\bf n},k,r} \int x^{\nu +r} \widehat{s}_{1,j}^2(x) \widehat{s}_{1,k}^1(x)d\sigma_0^2(x) =
\]
\[  \int x^{\nu} \widehat{s}_{1,j}^2(x) \sum_{k=0}^{m_2} a_{{\bf n},k}(x) \widehat{s}_{1,k}^1(x) d\sigma_0^2(x) = 0,
\]
as we needed to prove.

From Theorem \ref{teo:1}, we know that ${\mathcal{A}}_{\bf n}$ has at most $|{\bf n}_1| -1 = |{\bf n}_2|$ zeros on
${\mathbb{C}} \setminus \mbox{\rm Co}(\supp \sigma_1^1)$. From (\ref{orto}) it follows that this form has at least $|{\bf n}_2|$ sign changes in the interior of $\mbox{\rm Co}(\supp \sigma_0^2)$. In fact, assume that ${\mathcal{A}}_{\bf n}$ has exactly $N \leq |{\bf n}_2| -1$ sign changes in the indicated set at the points $x_1,\ldots,x_N$. According to Theorem \ref{teo:1}, $(1,\widehat{s}^2_{1,1},\ldots,\widehat{s}^2_{1,m_2})$ is also an AT system. Using the properties of Tchebychev systems (see \cite[Theorem 1.3]{KN}, we can find polynomials $p_0,\ldots,p_{m_2},$  with $\deg p_j \leq n_{2,j} -1,$ such that ${\mathcal{L}}_{\bf n}(x) = p_0(x) + \sum_{j=1}^{m_2} p_{j}(x) \widehat{s}_{1,j}^2(x)$
changes  sign at $x_1,\ldots,x_N,$ and has no other points where it changes sign in the interior of $\mbox{Co}(\supp \sigma^2_0)$. Thus, the function ${\mathcal{L}}_{\bf n}(x) {\mathcal{A}}_{{\bf n}}(x)$ has constant sign on $\mbox{Co}(\supp \sigma^2_0)$. This contradicts (\ref{orto}) since $\sigma^2_0$ is a measure with constant sign whose support contains infinitely many points. So the number of sign changes is $|{\bf n}_2|$ as claimed.   \hfill $\Box$ \vspace{0,2cm}

{\bf Proof of Theorem \ref{teo:2}.}  Suppose that for some ${\bf n}, {\mathbb{A}}_{\bf n}$ is not normal. That is,  some component $a_{{\bf n},k}$ of ${\mathbb{A}}_{\bf n}$ has $\deg a_{{\bf n},k} \leq n_{1,k} -2$. Then, according to Theorem \ref{teo:1}, ${\mathcal{A}}_{\bf n}$ can have on the interval $\mbox{Co}(\supp \sigma^2_0)$ at most $|{\bf n}_1| -2 =  |{\bf n}_2| -1$ zeros, but we know from Theorem \ref{lem:orto} that this is not the case, so ${\mathbb{A}}_{\bf n}$ must be normal and perfectness has been established.

Let us assume that there are two non collinear solutions ${\mathbb{A}}_{\bf n}, {\mathbb{A}}_{\bf n}^*,$ to a)-b).   Then, there exists a real constant $C \neq 0$ such that ${\mathbb{A}}_{\bf n} - C {\mathbb{A}}_{\bf n}^* \not \equiv 0$ and at least one of the components  of ${\mathbb{A}}_{\bf n} - C {\mathbb{A}}_{\bf n}^* $ satisfies $\deg (a_{{\bf n},k } - C a_{{\bf n},k }^*) \leq n_{1,k} -2.$ This is not possible since ${\mathbb{A}}_{\bf n} - C{\mathbb{A}}_{\bf n}^*$ also solves a)-b) and according to what was proved above all its components must have maximum possible degree. \fp \vspace{0,2cm}

\begin{rem} \label{moment} In Theorems \ref{teo:1}-\ref{teo:2} we have imposed the smallest number of moments on the generating measures which guarantees the results (using the method we have devised) regardless the multi-index under consideration. For example, if the multi-index is strictly increasing then that amount is necessary to completely carry out the inversion process on the measures. But, as we saw, for example in Lemma \ref{lem:3a}, for specific multi-indices a smaller amount of moments may be sufficient. In any case, this is not a matter of utmost importance since in applications the usual situation is that the measures have all their moments finite.
\end{rem}

{\bf Proof of Corollary \ref{biort}}. Given ${\bf n}_1 \in \Lambda_1$, take ${\mathcal{Q}}_{{\bf n}_1} = {\mathcal{A}}_{\bf n}$ where ${\bf n} = ({\bf n}_1;{\bf n}_2)$ and ${\bf n}_2$ is the unique element in $\Lambda_2$ such that $|{\bf n}_1| = |{\bf n}_2| +1$. Interchange the roles of $S_1$ and $S_2$ and for each ${\bf n}_2 \in \Lambda_2$ take ${\mathcal{P}}_{{\bf n}_2} = {\mathcal{A}}_{\bf n}$ where ${\bf n} = ({\bf n}_2;{\bf n}_1)$ and ${\bf n}_1$ is the unique element in $\Lambda_1$ such that $|{\bf n}_2| = |{\bf n}_1| +1$. Since the sequences of multi-indices $\Lambda_1$ and $\Lambda_2$ are ordered (\ref{bi1}) is a consequence of (\ref{orto}). The statement about the location of the zeros of the bi-orthogonal forms is contained in the last statement of Theorem \ref{lem:orto}. The assertions concerning the degrees of the polynomials $q_{{\bf n}_1,k}, p_{{\bf n}_2,k}$, the uniqueness of the forms except for a constant factor, and (\ref{bi2}) follow from the perfectness proved in Theorem \ref{teo:2}. \hfill $\Box$ \vspace{0,2cm}

A repeated use of Lemma \ref{lem:4} allows to give that result a more conclusive form. The proof is exactly the same as that of \cite[Theorem 1.3]{LF4} so we limit ourselves to the statement.

\begin{teo}
\label{teo:3}
Let $ (s_{1,1},\ldots,s_{1,m}) = {\mathcal{N}}(\sigma_{1},\ldots,\sigma_m)$ be such that the first $2^k +1$ moments of $\sigma_k,k=1,\ldots,m$ are finite, and ${\bf n} = (n_0,\ldots,n_m) \in {\mathbb{Z}}_+^{m+1}$. Then, there exists a permutation $\lambda$ of $(0,\ldots,m)$ which reorders the components of ${\bf n}$ decreasingly, $n_{\lambda(0)} \geq \cdots \geq n_{\lambda(m)},$ and an associated Nikishin system $S(\lambda) = (r_{1,1},\ldots,r_{1,m}) = {\mathcal{N}}(\rho_{1},\ldots,\rho_m)$ such that for any real polynomials
$p_k, \deg p_k \leq n_k -1$, there exist real polynomials $q_k$ such that
\[ p_0 + \sum_{k=1}^m p_k \widehat{s}_{1,k} = (q_0 + \sum_{k=1}^m q_k \widehat{r}_{1,k})\widehat{s}_{1,\lambda(0)}, \qquad \deg q_k \leq n_{\lambda(k)} -1, \qquad k=0,\ldots,m.
\]
\end{teo}

Here $\widehat{s}_{1,0}$ denotes the function identically equal to $1$; this is relevant when $\lambda(0) =0$. We do not know if there is an $S(\lambda)$ for each $\lambda$ which reorders the components of ${\bf n}$ decreasingly. We can say that there exists $S(\lambda)$ (but not exclusively) for that $\lambda$ which  additionally verifies that for all $0 \leq j <k \leq n$ with $n_j = n_k$ then also $\lambda(j) < \lambda(k)$. \vspace{0,2cm}

Likewise, following the proof of \cite[Theorem 1.4]{LF4} one obtains

\begin{teo}\label{teo:4} Let $S^1 =   {\mathcal{N}}(\sigma_0^1,\ldots,\sigma_{m_1}^1),  S^2= {\mathcal{N}}(\sigma_0^2,\ldots,\sigma_{m_2}^2),$ be two compatible Nikishin systems such that the first $2^k +1$  moments of $\sigma_k^i, k=1,\ldots,m_i, i=1,2,$ are finite and all the moments of $\sigma_0^2 (= \sigma_0^1)$ are finite. Let ${\bf n} \in {\mathbb{Z}}_+^{m_1 +1}\times {\mathbb{Z}}_+^{m_2 +1}, |{\bf n}_1| = |{\bf n}_2| +1,$ be given. Denote by $\lambda_2$ and $S(\lambda_2)= {\mathcal{N}}(\rho_1^2,\ldots,\rho_{m_2}^2)$ a permutation and a Nikishin system associated with ${\mathcal{N}}(\sigma_1^2,\ldots,\sigma_{m_2}^2)$ and ${\bf n}_2$ by Theorem $\ref{teo:3}$ with the reordering effect.  Construct $(r^2_{0,0},\ldots,r^2_{0,m_2})= {\mathcal{N}}(\rho_0^2,\ldots,\rho_{m_2}^2)$, where  $\rho^2_0 = \widehat{s}^2_{1,\lambda_2 (0)}\sigma_0^2$.  Then
\[ \int x^{\nu} {\mathcal{A}}_{{\bf n}}(x) d r_{0,k}^2(x) =0,\qquad \nu =0,\ldots,n_{2,\lambda_2 (k)} -1, \qquad k=0,\ldots,m_2.
\]
\end{teo}

Let $l = (l_1,l_2), 0 \leq l_1 \leq m_1, 0 \leq l_2 \leq m_2,$ be a pair of integers and ${\bf n} = ({\bf n}_1;{\bf n}_2) \in {\mathbb{Z}}_+^{m_1 +1}\times {\mathbb{Z}}_+^{m_2 +1}, |{\bf n}_1| = |{\bf n}_2| +1,$ be given. By ${\bf n}^l = ({\bf n}_1^{l_1};{\bf n}_2^{l_2})$ we denote the multi-index obtained adding one to the $l_i +1$ component of ${\bf n}_i, i=1,2$.

The zeros of "consecutive" linear forms ${\mathcal{A}}_{\bf n}$ interlace. Theorems \ref{teo:3} and \ref{teo:4} allow to reduce the proof to the case when the components of ${\bf n}_1$ and ${\bf n}_2$ are decreasing. In that situation, the proof is similar to that of \cite[Theorem 3.5]{FLLS}) for bounded and non-intersecting supports in the generating measures.

\begin{cor} \label{cor:2}  Let $S^1 =   {\mathcal{N}}(\sigma_0^1,\ldots,\sigma_{m_1}^1),  S^2= {\mathcal{N}}(\sigma_0^2,\ldots,\sigma_{m_2}^2),$ be two compatible Nikishin systems such that the first $2^k +1$  moments of $\sigma_k^i, k=1,\ldots,m_i, i=1,2,$ are finite and all the moments of $\sigma_0^2 (= \sigma_0^1)$ are finite. Let ${\bf n} \in {\mathbb{Z}}_+^{m_1 +1}\times {\mathbb{Z}}_+^{m_2 +1}, |{\bf n}_1| = |{\bf n}_2| +1.$  By $\lambda_i$ and $S(\lambda_i)= {\mathcal{N}}(\rho_1^i,\ldots,\rho_{m_i}^i), i=1,2,$ denote permutations and Nikishin systems associated with ${\mathcal{N}}(\sigma_1^i,\ldots,\sigma_{m_i}^i)$ and ${\bf n}_i$ through Theorem $\ref{teo:3}$, respectively, with the reordering effect. Let $l$ be such that the same permutations and systems work  if we replace ${\bf n}_i$ by ${\bf n}_i^{l_i}$. Then, between two consecutive zeros of ${\mathcal{A}}_{{\bf n}^l}$ in the interior of $\mbox{\rm Co}(\supp \sigma_0^1)$ lies exactly one zero of ${\mathcal{A}}_{\bf n}$.
\end{cor}

\begin{rem} It is easy to construct complete sequences of multi-indices $\Lambda_1, \Lambda_2,$ for which Corollary \ref{cor:2} is applicable to any two consecutive multi-indices in $\Lambda_1$ and $\Lambda_2$. In this situation the zeros of any two consecutive forms ${\mathcal{Q}}_{{\bf n}_1}$ and  ${\mathcal{P}}_{{\bf n}_2}$ interlace.
\end{rem}

Let us restrict our attention to type II approximation and the proof of Corollary \ref{cor:1}.

\begin{defi}
\label{def:hausdroff} Let $E$ be a subset of the complex plane and $\mathcal{U}$ the class of all coverings of $E$ by disks $U_n$. The radius of $U_n$ is denoted $|U_n|$. The (one dimensional) Hausdorff content of $E$ is
\[ h(E) = \inf \{\sum|U_n|: \{U_n\} \in {\mathcal{U}}\}.
\]
\end{defi}

Let $\{f_n\}_{n \in \Lambda}$ be a sequence of functions defined on a region $D \subset {\mathbb{C}}$. We say that $\{f_n\}_{n \in \Lambda}$ converges to $f$ in Hausdorff content on $D$ if for every compact set ${\mathcal{K}}  \subset D$ and any $\varepsilon > 0$
\[ \lim_{n \in \Lambda} h(\{z \in {\mathcal{K}}: |f_n(z) - f(z)| > \varepsilon\}) =0.
\]
We denote this by
\[ {\mathcal{H}}-\lim_{n \to \infty} f_n = f, \qquad {\mathcal{K}} \subset D.
\]

In \cite[Lemma 1]{Gon}, A.A. Gonchar proved that if the functions $f_n$ are holomorphic in $D$ and they converge in Hausdorff content to $f$ in $D$, then $f$ is in fact holomorphic in $D$ (more precisely, differs from a holomorphic function on a set of zero Hausdorff content) and the convergence (to the equivalent holomorphic function) is  uniform on each compact subset of $D$.
\vspace{0.2cm}

{\bf Proof of Corollary \ref{cor:1}.} In \cite[Theorem 1]{LF3}  it was proved (see also Remark 3 at the end of that paper) that under the assumptions of the corollary, for each $k=0,\ldots,m,$
\[ {\mathcal{H}}-\lim_{{\bf n} \in \Lambda} R_{{\bf n},k} = \widehat{s}_k, \qquad {\mathcal{K}} \subset \overline{\mathbb{C}} \setminus \mbox{Co}(\supp \sigma_0).
\]
Due to Gonchar's lemma and the last assertion of Theorem \ref{lem:orto}, it follows that convergence is uniform on each compact subset of $\overline{\mathbb{C}} \setminus \mbox{Co} (\supp \sigma_0)$. Regarding the proof of the rate of convergence, we refer to \cite[Corollary 1]{LF3} and the last sentence on page 104 of  the same paper.  \hfill $\Box$

\end{document}